\documentclass[a4paper,11pt]{article}
\usepackage[active]{srcltx}
\usepackage{amsmath}\usepackage{amssymb}\usepackage{amscd}

\newcommand{\pos}{{\mathrm{pos}}}

\newcommand{\Tr}{\textrm{Tr}}

\newcounter{pic}
\newenvironment{pic}[1][\bf Fig. \arabic{pic}]{
        
        \refstepcounter{pic}\noindent\textbf{#1.}${}$\hspace{5pt}${}$\it}{}
\newtheorem{defi}{Definition}\newtheorem{prop}[defi]{Proposition}\newtheorem{lemm}[defi]{Lemma}      
\newtheorem{coro}[defi]{Corollary}\newtheorem{rema}[defi]{Remark}

\newcommand{\E}{\mathfrak{E}}
\newcommand{\blambda}{\boldsymbol{\lambda}}
\newcommand{\balpha}{\boldsymbol{\alpha}}
\newcommand{\wbg}{w_{\blambda}^{\gamma}}
\newcommand{\wbgo}{w_{\blambda}^{\gamma^{\circ}}}
\newcommand{\twbg}{\widetilde{w}_{\blambda}^{\gamma}}

\newcommand{\B}{\mathcal{B}}

\usepackage{empheq}

\begin{document}

$\ $

\begin{center}
\bigskip\bigskip

{\Large\bf 
Induced representations and traces for chains of affine and cyclotomic Hecke algebras} 

\vspace{1.6cm}
{\large {\bf O. V. Ogievetsky$^{\circ
}$\footnote{On leave of absence from P. N. Lebedev Physical Institute, Leninsky Pr. 53,
117924 Moscow, Russia}
and L. Poulain d'Andecy$^{\diamond
}$\footnote{The second author was supported by ERC-advanced grant no. 268105.}}}

\vskip 1cm
$\circ\ ${Center of Theoretical Physics\\
Aix Marseille Universit\'e, CNRS, UMR 7332, 13288 Marseille, France\\
Universit\'e de Toulon, CNRS, UMR 7332, 83957 La Garde, France}

\vskip .4cm
$\diamond\ ${Korteweg-de Vries Institute for Mathematics, University of Amsterdam\\
P.O. Box 94248, 1090 GE Amsterdam, The Netherlands}
\end{center}

\vskip 1cm
\begin{abstract} \noindent
Properties of relative traces and symmetrizing forms on chains of cyclotomic and affine Hecke algebras are studied. The study relies on a use of 
bases of these algebras which generalize a normal form for elements of the complex reflection groups $G(m,1,n)$, $m=1,2,\dots,\infty$, constructed by a recursive use of the Coxeter--Todd algorithm. Formulas for inducing, from representations of an algebra in the chain, representations of the next member of the chain are presented.  
\end{abstract}

\vskip 1.2cm
\section{{\hspace{-0.55cm}.\hspace{0.55cm}}Introduction\vspace{.25cm}}

The Coxeter--Todd algorithm \cite{CT} is a powerful tool for constructing a normal form for group elements with respect to a given subgroup. 
For an ascending chain of groups the Coxeter--Todd algorithm establishes recursively a global normal form for group elements. We apply the 
algorithm to the chain, in $n$, of the complex reflection groups $G(m,1,n)$. The algorithm works uniformly for all $n>1$. 
The normal form of an element of $G(m,1,n)$ is $\textsc{l}w$, where $w\in G(m,1,n-1)$  
and $\textsc{l}$ runs through a left transversal of the subgroup $G(m,1,n-1)$ in the group $G(m,1,n)$; moreover, the left transversal is uniform as well.  
In this sense (that the normal form is $\textsc{l}w$), this normal form is ``inductive'', in contrast to the normal forms used in \cite{AK,Bre-M}; 
as in \cite{Bre-M}, this normal form consists of reduced expressions in terms of generators of $G(m,1,n)$. 

\vskip .1cm
We allow the value $m=\infty$. The group $G(\infty ,1,n)$ is the affine Weyl group of type GL. 
The Coxeter--Todd algorithm is applicable for $m=\infty$ and we find the corresponding left transversal and the normal form for group elements.

\vskip .1cm 
We denote by $H(m,1,n)$ the standard deformation of the complex reflection group $G(m,1,n)$; it has been introduced 
in \cite{AK,Bro-M,C2} and is called the Ariki--Koike algebra, or the cyclotomic Hecke algebra. Again, $m=\infty$ is allowed, it corresponds to the affine Hecke algebra $\hat{H}_n$ of type GL. 
The above inductive normal form has a nice generalization to a basis $\mathcal{B}$ of the algebra $H(m,1,n)$, $m=1,2,\dots ,\infty$, inductive with respect to the chain, in $n$. 
The use of the basis $\mathcal{B}$ allows for a representation-free approach, avoiding a dimension count, to the structure theory of the cyclotomic and affine Hecke algebras. 
The fact, that $\mathcal{B}$ is a basis, implies that $H(m,1,n)$ is a flat deformation of $\mathbb{C}G(m,1,n)$. 
The basis $\mathcal{B}$ is different from the bases in \cite{AK,Bre-M} (introduced there for $m<\infty$) but is similar to the basis in \cite{L} used for the study of Markov traces on the (cyclotomic) Hecke algebra. 
Nevertheless the results of \cite{L} rely on the basis in \cite{AK}, obtained as an outcome of the classification of irreducible representations of the algebra $H(m,1,n)$, for $m<\infty$.  
Our arguments do not refer to the representation theory and work for $m=\infty$ as well; 
the proof is done more in a spirit of classical proofs for the usual Hecke algebra
(note that a basis of cyclotomic quotients of the degenerate affine Hecke algebras is established in \cite{Kl} without the use of the representation theory). 

\vskip .1cm
A key ingredient in our study consists in explicit formulas for inducing representations of the algebra $H(m,1,n-1)$ to representations of the algebra $H(m,1,n)$, 
$m=1,2,\dots ,\infty$. A one-dimensional representation of the algebra $H(m,1,n-1)$ induces a natural analogue,  
for the algebra $H(m,1,n)$, of the Burau representation.

\vskip .1cm
The cyclotomic/affine Hecke algebras possess commutative Bethe subalgebras (see \cite{MTV} for the symmetric group case). The Bethe subalgebras play a fundamental role 
in the theory of integrable systems such as chain and Gaudin models, see, e.g., \cite{IK,IO}. One of ways to construct the Bethe subalgebras is based 
on ``relative traces", the linear maps $\Tr_k\colon H(m,1,k)\to H(m,1,k-1)$, satisfying certain conditions, see Section \ref{ConditExpec} for precise definitions. We use the basis 
$\mathcal{B}$ to establish the existence and uniqueness properties of the relative traces for the chain of the cyclotomic/affine Hecke algebras. 

\vskip .1cm 
A composition $\Tr:=\Tr_1\circ\dots\circ\Tr_{n-1}\circ\Tr_n$ of relative traces defines a Markov trace on the algebra $H(m,1,n)$ and all Markov traces constructed in \cite{L} can be thus obtained.
In particular, the existence and unicity of the Markov traces constructed in \cite{L} is reobtained in a different way via the study of the relative traces.
If $\Tr$ is non-degenerate on both $H(m,1,n)$ and $H(m,1,n-1)$ then $\Tr_n$ coincides with the conditional expectation.

\vskip .1cm
 Besides, there is a family $L^\gamma$ of central forms on $H(m,1,n)$ which possess a 
multiplicativity property with respect to the basis $\mathcal{B}$. 
Up to the standard involution of $H(m,1,n)$, these central forms $L^{\gamma}$ are 
Markov traces on $H(m,1,n)$ with the ``Markov parameter" equal to 0 (see Section \ref{sec-symm} for precisions).
As an application of use of the basis $\mathcal{B}$, we obtain 
expressions for the products of the generators and the basis elements and give an independent direct proof of the centrality of 
the linear forms $L^\gamma$ on $H(m,1,n)$. 
{}For $m<\infty$, one of these central forms, denoted $L^{\gamma^{\circ}}$, coincides with the central form introduced in \cite{Bre-M} and further studied in \cite{MM}. In the terminology of \cite{MM}, the inductive basis $\mathcal{B}$ is \emph{quasi-symmetric} with respect to the central form $L^{\gamma^{\circ}}$, as are the different bases in \cite{AK,Bre-M}. 
We also note that, for each central form $L^{\gamma}$, an easy modification of the inductive basis $\mathcal{B}$ yields a basis quasi-symmetric with respect to $L^{\gamma}$.

\vskip .1cm 
In \cite{OPdA4}, a \emph{fusion formula} for a complete set of primitive idempotents of the algebra $H(m,1,n)$, for finite $m$, is obtained. The fusion formula is well adapted to the inductive basis $\mathcal{B}$ (for more details, see also \cite{OPdA3a}). For $m<\infty$, it turns out that, due to the multiplicativity property, the weights of the central forms $L^{\gamma}$ can be easily calculated with the help of the fusion formula for the algebra $H(m,1,n)$.
The weights of any Markov trace from \cite{L} have been calculated in \cite{GIM} (the weights of the form $L^{\gamma^{\circ}}$ have been calculated independently in \cite{M}, see also \cite{CJ}). The formulas for the weights obtained here are different from the formulas in \cite{GIM} and generalize, for any $L^{\gamma}$, the so-called ``cancellation-free" formula obtained in \cite{CJ} for the weights of $L^{\gamma^{\circ}}$.

\vskip .3cm
The paper is organized as follows.

\vskip .05cm 
In Section \ref{sec-group}, we present the Coxeter--Todd algorithm for the chain (with respect to $n$) of the groups $G(m,1,n)$. We establish the resulting normal form for elements of $G(m,1,n)$. 

\vskip .05cm
The normal form for elements of $G(m,1,n)$ suggests a basis for the algebra $H(m,1,n)$. In Sections \ref{sec-nf} and \ref{sec-nf2}, we show that this is 
indeed a basis. Several known facts about the chain (with respect to $n$) of the algebras $H(m,1,n)$ are reestablished with the help of this basis. In particular,
we show that $H(m,1,n)$ is a flat deformation of the group ring of $G(m,1,n)$. We also give the formulas for the induced representations. 

\vskip .05cm
Section \ref{ConditExpec} is devoted to the relative traces for the chain of the cyclotomic/affine Hecke algebras. 

\vskip .05cm
In Section \ref{sec-symm}, we complete the study of the multiplication of the generators on the basis elements. We then use the results to establish properties of central forms $L^\gamma$ on $H(m,1,n)$. For the cyclotomic Hecke algebras, 
we present a calculation, based on the fusion formula, 
of weights of the forms $L^\gamma$ in Subsection \ref{weightsfusion}.

\vskip .3cm\noindent
{\bf Notation.} 

\vskip .1cm
\noindent
$\E_m:=\{ 0,1,\dots ,m-1\}$ for finite $m$ and $\E_{\infty}:={\mathbb{Z}}$\,;

\vskip .1cm
\noindent
${\mathcal{A}}_m:=\mathbb{C}[q^{\pm 1},v_1^{\pm 1},\dots ,v_m^{\pm 1}]$ for finite 
$m$ and ${\mathcal{A}}_{\infty}:=\mathbb{C}[q,q^{-1}]$; here $q,v_1,\dots,v_m$ are indeterminates.

\vskip .1cm\noindent
The symbol $\square$ stands for an end of a proof, $\triangle$ - for an end of a remark.

\section{\hspace{-0.55cm}.\hspace{0.55cm}Normal form for the group $G(m,1,n)$}\label{sec-group}

Let $m\in \mathbb{Z}_{>0}\cup\{\infty\}$. The group $G(m,1,n)$ is generated by the elements $t$, $s_1$, $\dots$, $s_{n-1}$ with the defining relations:
\begin{equation}\label{def2a}
\left\{\begin{array}{ll}
s_i^2=1 & \textrm{for $i=1,\dots,n-1$\ ,}\\[.2em]
s_is_{i+1}s_i=s_{i+1}s_is_{i+1} & \textrm{for $i=1,\dots,n-2$\ ,}\\[.2em]
s_is_j=s_js_i & \textrm{for all $i,j=1,\dots,n-1$ such that $|i-j|>1$\ ,}
 \end{array}\right.
\end{equation}
and
\begin{equation}\label{def2b}
\left\{\begin{array}{ll}
t^m=1 & \textrm{if $m<\infty$\ ,}\\[.2em]
ts_1ts_1=s_1ts_1t\ ,\\[.2em]
ts_i=s_it & \textrm{for $i=2,\dots,n-1$\ .}
\end{array}\right.
\end{equation}
Let $t_1:=t$ and, inductively, $t_{i+1}:=s_it_is_i$ for $i=1,\dots,n-1$. The following holds (see, e.g., \cite{OPdA3}) 
\begin{equation}\label{rel-wr} 
\left\{\begin{array}{ll}
t_a^m=1 & \textrm{for $a=1,\dots,n$,\ \ if\ \ $m<\infty$\ ,}\\[.2em]
t_at_b=t_bt_a & \textrm{for $a,b=1,\dots,n$\ ,}\\[.2em]
s_it_a=t_{\pi_i(a)}s_i \ \ & \textrm{for $i=1,\dots,n-1\,,\ a=1,\dots,n$\ ,}
\end{array}\right.
\end{equation}
where $\pi_i$ is the transposition $(i,i+1)$. The relations (\ref{rel-wr}), together with (\ref{def2a}), are the defining relations of the wreath product $C_m\wr S_n$ of the cyclic group $C_m$ with $m$ elements (the infinite cyclic group if $m=\infty$) by the symmetric group $S_n$. The relations (\ref{def2b}) are included in the relations (\ref{rel-wr}).
Thus the group $G(m,1,n)$ is isomorphic to the group $C_m\wr S_n$. Let $\gamma$ be a generator of the group $C_m$. The 
images of the generators $t$, $s_1$, $\dots$, $s_{n-1}$ of the group $G(m,1,n)$ in $C_m\wr S_n$ are 
\begin{equation}\label{sta-iso}
t\mapsto \Bigl(\left(\begin{array}{c}\gamma\\ e\\ \vdots\\ e\end{array}\right),e\Bigr)\ ,
\quad s_i\mapsto \Bigl(\left(\begin{array}{c}e\\ e\\ \vdots\\ e\end{array}\right),(i,i+1)\Bigr)\ ,
\end{equation}
where $e$ is the unit element in a corresponding subgroup and the vectors are elements of the Cartesian product of $n$ copies of $C_m$.
Besides, the subgroup generated by the elements $t$, $s_1$, $\dots$, $s_{n-2}$ is isomorphic to $G(m,1,n-1)$ and so we have the chain, in $n$, 
of groups $G(m,1,n)$. 

\begin{rema}{\hspace{-.2cm}.\hspace{.02cm}}
\rm{
We recall that
the group $G(m,1,n)$, $m=1,\dots,\infty$, admits the following description. Let $E$ be the set $\{(\varrho ,a)\ \vert\ \varrho\in C_m,\,\,a=1,\dots,n\}$. Define the 
action of $C_m$ on this set: $\varrho\cdot(\varrho',a)=(\varrho\varrho',a)$, $\varrho,\varrho'\in C_m,\,\,a=1,\dots,n$. Denote by $\text{\sf Perm}(E)$ the group of permutations of the set $E$. Let $\text{\sf Perm}^0(E)$ be the subgroup of $\text{\sf Perm}(E)$ consisting of elements $\pi\in \text{\sf Perm}(E)$ such that:
\begin{equation}\pi(\varrho ,a)=\varrho\cdot\pi(e,a)\quad\textrm{for $\varrho\in C_m$ and $a=1,\dots,n$\ .}\end{equation}
To specify an 
element $\pi$ of $\text{\sf Perm}^0(E)$, it is enough to give the images under $\pi$ of the elements of 
the set $\{(e,a)\ \vert\ a=1,\dots,n\}$. 

\vskip .2cm
The group $G(m,1,n)$ is isomorphic to $\text{\sf Perm}^0(E)$. Indeed, let $\phi$ be the map from 
the set of generators of $G(m,1,n)$ to $\text{\sf Perm}^0(E)$ defined by:
\begin{equation}\begin{array}{ll}\phi(t)(e,a)=\left\{\begin{array}{ll}(\gamma,a) & \textrm{if $a=1$\ ,}\\[.2em]
(e,a) & \textrm{otherwise\, ;}\end{array}\right.  & \\ [1.5em]
\phi(s_i)(e,a)=(e,\pi_i(a))\quad\textrm{for $a=1,\dots,n$ and $i=1,\dots,n-1$\ }\\
                                              \end{array}\end{equation}
(recall that $\pi_i$ is the transposition $(i,i+1)$). 
The images of the elements $t,s_1,\dots ,s_{n-1}$, the permutations $\phi(t),\phi(s_1),\dots,\phi(s_{n-1})\in \text{\sf Perm}^0(E) $, satisfy the defining relations of 
the group $G(m,1,n)$, so the map $\phi$ extends to a homomorphism $G(m,1,n)\to\text{\sf Perm}^0(E)$. This homomorphism is surjective; it is a consequence of the following fact: for any $j=1,\dots,n$, we have
\begin{equation}
\phi (t_{j})(e,a)=\left\{\begin{array}{ll}(\gamma,a) & \textrm{if $a=j$\ ,}\\[.3em]
(e,a) & \textrm{otherwise\ .}\end{array}\right.\end{equation}
Using the isomorphism between $G(m,1,n)$ and $C_m\wr S_n$, we see that $\phi$ is injective: $\phi (\vec{v},\omega )(e,a)=
(v_{\omega (a)},\omega (a))$ where $(\vec{v},\omega )\in C_m\wr S_n$ and $a=1,\dots,n$, so the kernel of $\phi$ is trivial.
Therefore, the map $\phi$ is an isomorphism.
\hfill$\triangle$ 
}\end{rema}

\subsection{\hspace{-0.50cm}.\hspace{0.50cm}Coxeter--Todd algorithm for the chain of groups $G(m,1,n)$}

For a group $G$ given by generators and relations and a subgroup $W$ of $G$ generated by some subset of the generators of $G$, the Coxeter--Todd algorithm consists \cite{CT} in constructing the set of the left cosets of $W$ in $G$ and the action of the generators on this set. To every coset, a vertex in the 
Coxeter--Todd figure is associated; arrows stand for the action of the generators. The algorithm starts with the left coset $eW=W$ ($e$ is the identity element); 
only the generators of $G$ which are not in $W$ act non-trivially on this coset producing new vertices. At each step, using the relations of a given presentation we analyze the action of the generators on vertices which are already in the figure and add or identify 
possible cosets. The algorithm terminates when we know the action of all 
generators on every coset in the figure. The figure gives an upper bound for the order of a group $G$. The Coxeter--Todd algorithm lists the set of cosets and 
provides thereby a ``normal form of an element of $G$ with respect to $W$". 

\paragraph{1. Coxeter--Todd  figure for the chain.} Let $W$  be the subgroup of $G(m,1,n)$ generated by the elements $t,s_1,\dots,s_{n-2}$.
We present on Fig. \ref{CoxToddCyclo} the Coxeter--Todd figure for the group $G(m,1,n)$ with respect to its subgroup $W$.
\vskip .2cm
\begin{center}
 \setlength{\unitlength}{2500sp}
\begingroup\makeatletter\ifx\SetFigFontNFSS\undefined
\gdef\SetFigFontNFSS#1#2#3#4#5{  \reset@font\fontsize{#1}{#2pt}
  \fontfamily{#3}\fontseries{#4}\fontshape{#5}  \selectfont}
\fi\endgroup
\begin{picture}(12470,4456)(346,-4258)
{\thinlines
\put(541,-3346){\circle*{144}}}{\put(1486,-3346){\circle*{144}}}{\put(2476,-3346){\circle*{144}}}{\put(4609,-3346){\circle*{144}}}
{\put(5581,-3346){\circle*{144}}}{\put(6139,-2491){\circle*{144}}}{\put(5311,-2041){\circle*{144}}}{\put(12736,-3391){\circle*{144}}}
{\put(11746,-3391){\circle*{144}}}{\put(10801,-3391){\circle*{144}}}{\put(8614,-3391){\circle*{144}}}{\put(7669,-3391){\circle*{144}}}
{\put(3376,-871){\circle*{144}}}{\put(2539,-376){\circle*{144}}}{\put(1774,119){\circle*{144}}}{\put(586,-3346){\line( 1, 0){810}}}
{\put(1576,-3346){\line( 1, 0){810}}}{\put(2566,-3346){\line( 1, 0){270}}}{\put(2971,-3346){\line( 1, 0){270}}}
{\put(3376,-3346){\line( 1, 0){270}}}{\put(3781,-3346){\line( 1, 0){270}}}{\put(4186,-3346){\line( 1, 0){270}}}{\put(4681,-3346){\line( 1, 0){810}}}
{\put(5382,-2067){\line( 5,-3){694.853}}}{\put(5815,-3859){\vector(-1, 2){225}}}{\put(11836,-3391){\line( 1, 0){810}}}{\put(10891,-3391){\line( 1, 0){810}}}
{\put(10486,-3391){\line( 1, 0){270}}}{\put(10081,-3391){\line( 1, 0){270}}}{\put(9586,-3391){\line( 1, 0){270}}}{\put(9181,-3391){\line( 1, 0){270}}}
{\put(8686,-3391){\line( 1, 0){270}}}{\put(7741,-3391){\line( 1, 0){810}}}{\put(6211,-2491){\line( 1, 0){270}}}{\put(5918,-3953){\line( 1,-1){270}}}
{\put(6301,-4246){\line( 1, 0){450}}}{\put(7642,-3409){\line(-1,-2){270}}}{\put(6886,-4201){\line( 2, 1){360}}}{\put(7408,-2905){\vector( 1,-2){225}}}
{\put(6661,-2491){\line( 1, 0){315}}}{\put(2630,-417){\line( 5,-3){694.853}}}{\put(1826, 42){\line( 5,-3){694.853}}}{\put(5546,-3353){\vector( 2, 3){540}}}
{\put(3421,-916){\line( 5,-3){450}}}{\put(4006,-1276){\line( 5,-3){450}}}{\put(4681,-1681){\line( 5,-3){450}}}{\put(7066,-2536){\line( 1,-1){270}}}
\put(61,-3706){\makebox(0,0)[lb]{\smash{{\SetFigFontNFSS{12}{14.4}{\rmdefault}{\mddefault}{\updefault}{$W$}}}}}
\put(440,-3256){\makebox(0,0)[lb]{\smash{{\SetFigFontNFSS{12}{14.4}{\rmdefault}{\mddefault}{\updefault}{$s_{n-1}$}}}}}
\put(1440,-3256){\makebox(0,0)[lb]{\smash{{\SetFigFontNFSS{12}{14.4}{\rmdefault}{\mddefault}{\updefault}{$s_{n-2}$}}}}}
\put(4680,-3256){\makebox(0,0)[lb]{\smash{{\SetFigFontNFSS{12}{14.4}{\rmdefault}{\mddefault}{\updefault}{$s_1$}}}}}
\put(5480,-2221){\makebox(0,0)[lb]{\smash{{\SetFigFontNFSS{12}{14.4}{\rmdefault}{\mddefault}{\updefault}{$s_1$}}}}}
\put(5630,-3076){\makebox(0,0)[lb]{\smash{{\SetFigFontNFSS{12}{14.4}{\rmdefault}{\mddefault}{\updefault}{$t$}}}}}
\put(5506,-3661){\makebox(0,0)[lb]{\smash{{\SetFigFontNFSS{12}{14.4}{\rmdefault}{\mddefault}{\updefault}{$t$}}}}}
\put(11700,-3256){\makebox(0,0)[lb]{\smash{{\SetFigFontNFSS{12}{14.4}{\rmdefault}{\mddefault}{\updefault}{$s_{n-1}$}}}}}
\put(10680,-3256){\makebox(0,0)[lb]{\smash{{\SetFigFontNFSS{12}{14.4}{\rmdefault}{\mddefault}{\updefault}{$s_{n-2}$}}}}}
\put(7730,-3256){\makebox(0,0)[lb]{\smash{{\SetFigFontNFSS{12}{14.4}{\rmdefault}{\mddefault}{\updefault}{$s_1$}}}}}
\put(1800,-61){\makebox(0,0)[lb]{\smash{{\SetFigFontNFSS{12}{14.4}{\rmdefault}{\mddefault}{\updefault}{$s_{n-1}$}}}}}
\put(2690,-556){\makebox(0,0)[lb]{\smash{{\SetFigFontNFSS{12}{14.4}{\rmdefault}{\mddefault}{\updefault}{$s_{n-2}$}}}}}
\put(7000,-3211){\makebox(0,0)[lb]{\smash{{\SetFigFontNFSS{12}{14.4}{\rmdefault}{\mddefault}{\updefault}{$t$}}}}}
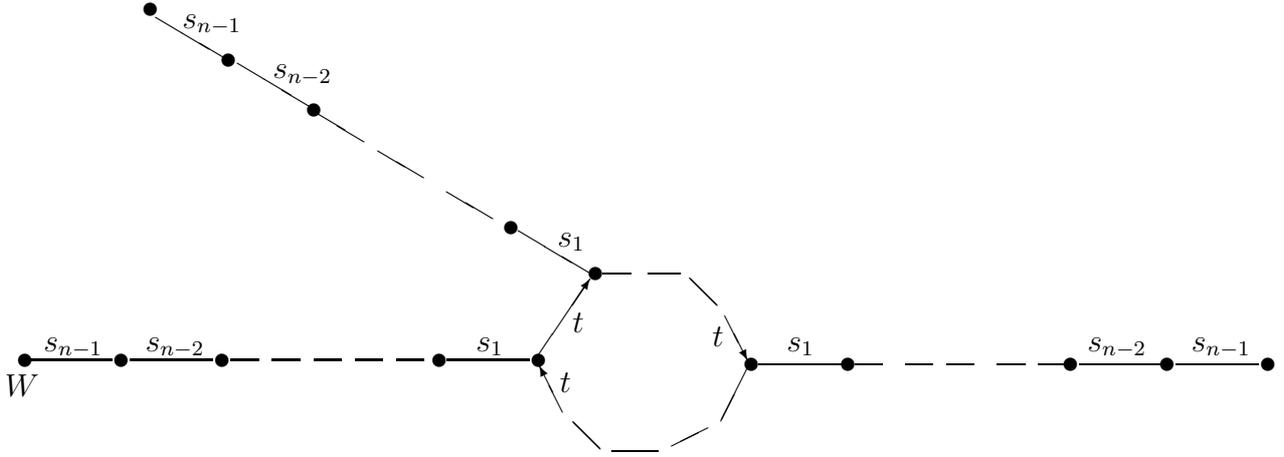
\put(3100,-5000){\begin{pic}\label{CoxToddCyclo}Coxeter--Todd figure for $G(m,1,n)$ with respect to $W$
\end{pic}}
\end{picture}
\end{center}

\vspace{1.2cm}
The action is indicated by oriented edges (labelled by the generators). For a generator of order 2 $\bigl($these are the generators 
$s_1,\dots,s_{n-1}$ of $G(m,1,n)$ if $m\neq 2$ $\bigr)$ 
an unoriented edge represents a pair of edges with opposite orientations. If a generator leaves a coset invariant, the corresponding edge is a 
loop which starts and ends at the vertex representing the coset. For brevity, these loops are omitted; only non-trivial actions are drawn.

\vskip .2cm
In the middle of Fig. 1 there is an $m$-gon with edges labelled by $t$ (if $m=\infty$, the $m$-gon is replaced by an infinite in two directions line with edges labelled by $t$ in one direction and by $t^{-1}$ in the other). At each vertex of the $m$-gon (respectively, of the infinite line if $m=\infty$) starts a tail with $n-1$ edges.
The tails from all vertices of the $m$-gon (respectively, of the infinite line) are labelled in the same way.
 
\subsection{\hspace{-0.50cm}.\hspace{0.50cm}Normal form for elements of $G(m,1,n)$.} 

The Coxeter--Todd algorithm leads to a normal form for the elements of $G(m,1,n)$ with respect to $G(m,1,n-1)$, $m=1,2,\dots ,\infty$.

\begin{prop}
{\hspace{-.2cm}.\hspace{.2cm}}\label{normalform-G} Any element $x\in G(m,1,n)$ can be written in the form 
\begin{equation}\label{normalform}x=s_js_{j-1}\dots s_1t^{\alpha}s_1s_2\dots s_{n-1}w\ ,
\end{equation}
where $j\in\{0,\dots,n-1\}$, $\alpha\in \E_m$ and $w\in W\simeq G(m,1,n-1)$ (by convention, the empty product is equal to $1$).
\end{prop}

{} For any $n$, given a certain normal form for elements of $G(m,1,n-1)$, the Proposition \ref{normalform-G} provides an induced normal form for elements of $G(m,1,n)$. 
Starting with a normal form (powers of $t$) for $G(m,1,1)\cong C_m$ we build, increasing $n$, 
recursively the global normal form for elements of $G(m,1,n)$ for any $n$.

\begin{rema}{\hspace{-.2cm}.\hspace{.2cm}}\label{reduced} {\rm The global normal form for elements of $G(m,1,n)$, obtained by the recursive application of the Proposition 
\ref{normalform-G}, provides, for finite $m$, 
expressions of minimal length for elements of $G(m,1,n)$ in terms of positive powers of the generators $t,s_1,\dots,s_{n-1}$. This can be proved 
using the same kind of argument as in \cite{Bre-M}. 
Or, one can directly use the results of \cite{Bre-M}. Namely, any element of $G(m,1,n)$ can be written in the form
\begin{equation}\label{nm-BM}\pi t_{n,a_n}\dots t_{2,a_2}t_{1,a_1},\ \quad\text{
$a_1,\dots,a_n\in\{ 0,\dots,m-1\} $ and $\pi\in S_n$,}\end{equation}
where $S_n$ is the symmetric group on $n$ letters, seen as the subgroup of $G(m,1,n)$ generated by $s_1,\dots,s_{n-1}$, and, for $k=1,\dots,n$,
\[t_{k,a}:=\left\{\begin{array}{ll}t^a s_1\dots s_{k-1} & \text{for $a>0$,}\\[0.5em]
 1 & \text{for $a=0$.}\end{array}\right.\]
Moreover, if $\pi\in S_n$ is given in a reduced form then the expression (\ref{nm-BM}) is minimal. 

\vskip .2cm
Now we prove by recurrence on $n$ that any expression of the form (\ref{normalform}) can be transformed into an expression of the form (\ref{nm-BM}) using only the braid relations in $G(m,1,n)$. 
It is trivial for $n=1$. Assume that $n>1$ and that an element $x\in G(m,1,n)$ is written in the normal form (\ref{normalform}). By the induction hypothesis, the expression 
$w\in G(m,1,n-1)$ can be transformed, using only the braid relations in $G(m,1,n-1)$, into
\[w=\pi t_{n-1,a_{n-1}}\dots t_{1,a_1}\ ,\]
for some $a_1,\dots,a_{n-1}\in\{0,\dots,m-1\}$ and a reduced expression $\pi\in S_{n-1}$. Thus, we have, using only the braid relations in $G(m,1,n)$,
\[x=s_{j}s_{j-1}\dots s_1\overline{\pi}t_{n,\alpha}t_{n-1,a_{n-1}}\dots t_{2,a_2}t_{1,a_1},\]
where $\overline{\pi}$ is the element of $S_n$ obtained from the reduced expression $\pi$ by replacing $s_i$ by $s_{i+1}$, for $i=1,\dots,n-2$. As $\overline{\pi}$ does not contain $s_1$, the expression $s_{j}s_{j-1}\dots s_1\overline{\pi}$ is a reduced expression of an element of $S_n$. 
The assertion is proved and
we conclude in particular that (\ref{normalform}) is a minimal expression of the element $x$.
\hfill$\triangle$ 
}\end{rema}

\section{\hspace{-0.55cm}.\hspace{0.55cm}Normal form for $H(m,1,n)$. Preparation}\label{sec-nf}

\subsection{\hspace{-0.50cm}.\hspace{0.50cm}Definitions}\label{sec-def} 

Let $m=1,2,\dots,\infty$ and $n>0$. The 
Hecke algebra $H(m,1,n)$ is the associative algebra over ${\mathcal{A}}_m$ generated by $\tau,\tau^{-1},\sigma_1,\dots,\sigma_{n-1}$ with the defining relations:
\begin{empheq}[]{alignat=4}
\label{rel1a}&\sigma_i\sigma_{i+1}\sigma_i=\sigma_{i+1}\sigma_i\sigma_{i+1} &\hspace{1cm}& \textrm{for all $i=1,\dots,n-2$\ ,}\\
\label{rel1b}&\sigma_i\sigma_j=\sigma_j\sigma_i && \textrm{for all $i,j=1,\dots,n-1$ such that $|i-j|>1$\ ,}\\
\label{rel1'a}&\tau\sigma_1\tau\sigma_1=\sigma_1\tau\sigma_1\tau\ ,&&\\
\label{rel1'b}&\tau\sigma_i=\sigma_i\tau && \textrm{for $i>1$\ ,}\\
\label{rel1''a} & \sigma_i^2=(q-q^{-1})\sigma_i+1 & &  \textrm{for all $i=1,\dots,n-1$\ ,}\\
\label{rel1''b} & (\tau-v_1)\dots(\tau-v_m)=0 & & \textrm{if $m<\infty$\ .}
\end{empheq}
Note that for finite $m$, $\tau^{-1}$ can be excluded from the set of generators since the relation (\ref{rel1''b}) allows 
to express $\tau^{-1}$ as a linear combination, with coefficients in ${\mathcal{A}}_m$, of positive powers of $\tau$. 

\vskip .2cm
For finite $m$, the algebra $H(m,1,n)$ is the cyclotomic Hecke algebra associated to the complex reflection group $G(m,1,n)$. The algebra $H(\infty,1,n)$ is the affine Hecke algebra $\hat{H}_n$ of type GL. 
The cyclotomic Hecke algebra $H(m,1,n)$ can be seen as the quotient of the affine Hecke algebra $\hat{H}_n$ by 
the relation (\ref{rel1''b}).

The braid group $B_n$ of type A (or simply the braid group) on $n$ strands is the group generated by the elements $\sigma_1$, $\dots$, $\sigma_{n-1}$ with the defining relations (\ref{rel1a})--(\ref{rel1b}), and the braid group $\alpha B_n$ of type B (sometimes called \emph{affine} braid group) is the group generated by the elements $\tau$, $\sigma_1$, $\dots$, $\sigma_{n-1}$ with the defining relations (\ref{rel1a})--(\ref{rel1'b}). The affine Hecke algebra $\hat{H}_n$ of type GL is the quotient of the group algebra of the braid group $\alpha B_n$ of type B by the quadratic relations (\ref{rel1''a}). 

Recall that the usual Hecke algebra $H_n$ of type A is the algebra over $\mathbb{C}[q,q^{-1}]$ generated by the elements $\sigma_1$, $\dots$, $\sigma_{n-1}$ with the relations (\ref{rel1a})--(\ref{rel1b}) and (\ref{rel1''a}). Thus, in particular, $H(1,1,n)$ is isomorphic to 
$\mathcal{A}_1\otimes H_n$ (the tensor product is over $\mathbb{C}[q,q^{-1}]$).
For $m=2$, the algebra $H(2,1,n)$ is isomorphic to the Hecke algebra of type B. 

\vskip .2cm
The algebra $H(m,1,n)$, $m=1,2,\dots,\infty$, is a deformation of the group algebra $\mathbb{C}G(m,1,n)$.
Indeed, the group algebra $\mathbb{C}G(m,1,n)$ is the specialization  of $H(m,1,n)$ for the following numerical values of the parameters:
\begin{equation}\label{specialization}
q\mapsto\pm1\,,\ \ \ \text{and, for finite $m$,}\ \ \ v_j\mapsto\xi_j\,,\ j=1,\dots,m,
\end{equation}
where $\{\xi_j\}_{j=1,\dots,m}$ is the set of all $m$-th roots of unity. 

\vskip .2cm
As $n$ varies, the algebras $H(m,1,n)$ form an ascending chain, in $n$, of algebras:
\begin{equation}\label{chaine}H(m,1,0):= \mathcal{A}_m\subset H(m,1,1)\subset\dots\subset H(m,1,n)\subset\dots\end{equation}
(it will be proved later that the elements $\tau,\tau^{-1}$ and $\sigma_1,\dots ,\sigma_{n-2}$ of the algebra $H(m,1,n)$ generate a subalgebra isomorphic to $H(m,1,n-1)$).
One has similar ascending chains of braid groups and affine braid groups.
Thus the reference to $n$ (as in $H_n$, $H(m,1,n)$, etc.) in the notation for the generators can be omitted.

\begin{rema}{\hspace{-0.2cm}.\hspace{0.2cm}} {\rm We use the same notation ``$\sigma_i$" for generators of different groups and algebras: these are braid and affine braid groups and 
Hecke, affine Hecke and cyclotomic Hecke algebras. The symbol $\tau$ is also used to denote a generator of several different objects.
This should not lead to any confusion, it will be clear from the context what is the algebra/group in question.
\hfill$\triangle$ 
}\end{rema}

\subsection{{\hspace{-0.50cm}.\hspace{0.50cm}}A spanning 
set of $H(m,1,n)$}

The normal form for elements of the algebra $H(m,1,n)$ we shall construct in several steps. This basis is an analogue 
of the normal form for the elements of the group $G(m,1,n)$ obtained in the previous Section. 
We start with the following Proposition. 

\begin{prop}
{\hspace{-.2cm}.\hspace{.2cm}}\label{normalform-H} Any $x\in H(m,1,n)$ can be written as a linear combination of elements
\begin{equation}\label{normalform2}
\sigma_j^{-1}\sigma_{j-1}^{-1}\dots \sigma_1^{-1}\tau^{\alpha}\sigma_1\sigma_2\dots \sigma_{n-1} w\ \ \text{with $j\in\{0,\dots,n-1\}$, $\alpha\in
\E_m$ and $w\in \tilde{W}$}
\end{equation}
where $\tilde{W}$ is the subalgebra generated by the elements $\tau,\tau^{-1},\sigma_1,\dots,\sigma_{n-2}$.
\end{prop}

\emph{Proof.} The unit element of the algebra $H(m,1,n)$ is in the set of elements (\ref{normalform2}), as well as all generators. We only have to check that 
the linear span of the elements (\ref{normalform2}) is stable under multiplication by the generators (then the linear span of the elements (\ref{normalform2}) is
a unital subalgebra, containing all generators, thus the whole algebra). We consider the left multiplication; we multiply an arbitrary element $\mathsf{E}$ of the 
form (\ref{normalform2}) by each generator $\mathsf{G}\in \{\sigma_1,\dots ,\sigma_{n-1},\tau,\tau^{-1}\}$ (recall that $\tau^{-1}$ is needed only if $m=\infty$) from the left and move in this product $\mathsf{GE}$ the 
original generator $\mathsf{G}$ to the right; the expression transforms; we follow the process until it becomes clear that the original product  
$\mathsf{GE}$ is a linear combination of elements of the form  (\ref{normalform2}).
 
\vskip .2cm
(i) We multiply the element (\ref{normalform2}) from the left by $\sigma_i$ with $i>j+1$. The element $\sigma_i$ commutes with 
$\sigma_j^{-1}\sigma_{j-1}^{-1}\dots \sigma_1^{-1}\tau^{\alpha}$ and thus moves through the combination 
$\sigma_j^{-1}\sigma_{j-1}^{-1}\dots \sigma_1^{-1}\tau^{\alpha}$ to the right without changes; then $\sigma_i$ moves through 
$\sigma_1\sigma_2\dots \sigma_{n-1}$ to the right, becoming $\sigma_{i-1}$; $\sigma_{i-1}w$ is again in $ \tilde{W}$ and we are done.

\vskip .2cm
(ii) We multiply the element (\ref{normalform2}) from the left by $\sigma_i$ with $i<j$. When $\sigma_i$ moves through 
$\sigma_j^{-1}\sigma_{j-1}^{-1}\dots \sigma_1^{-1}$, it transforms into $\sigma_{i+1}$; the element $\sigma_{i+1}$ commutes through $\tau^{\alpha}$ 
to the right without changes and then  $\sigma_{i+1}$ moves through $\sigma_1\sigma_2\dots \sigma_{n-1}$ to the right, 
becoming $\sigma_i$; as in (i), $\sigma_iw\in  \tilde{W}$.

\vskip .2cm
(iii) The assertion is immediate when we multiply the element (\ref{normalform2}) from the left by $\sigma_j$. 

\vskip .2cm
(iv) When we multiply the element (\ref{normalform2}) from the left by $\sigma_{j+1}$, for the proof it is enough to write $\sigma_{j+1}^{\phantom{-1}}
=\sigma_{j+1}^{-1}+(q-q^{-1})$.

\vskip .2cm
(v) We multiply the element (\ref{normalform2}) from the left by $\tau$. If $j=0$
one uses, if necessary, the characteristic equation (\ref{rel1''b})
for $\tau$ for finite $m$.
Let $j>0$. The element $\tau$ moves to the right until it reaches $\sigma_1^{-1}$. Then we use the Lemma \ref{lem3} below and obtain three terms. For the first term:
$$\begin{array}{ll}\sigma_j^{-1}\dots\sigma_2^{-1}\cdot\tau\sigma_1\tau^{\alpha}\cdot\sigma_2\dots \sigma_{n-1}\cdot w&= 
\sigma_j^{-1}\dots\sigma_2^{-1}\cdot\tau\cdot\sigma_1\sigma_2\dots \sigma_{n-1}\cdot\tau^{\alpha}w\\[.5em]
&=\tau\cdot\sigma_1\sigma_2\dots \sigma_{n-1}\cdot
\sigma_{j-1}^{-1}\dots\sigma_1^{-1}\tau^{\alpha}w\end{array}$$
and $\sigma_{j-1}^{-1}\dots\sigma_1^{-1}\tau^{\alpha}w\in\tilde{W}$. For the second term:
$$\sigma_j^{-1}\dots\sigma_2^{-1}\cdot\tau^{\alpha+1}\sigma_1\sigma_2\dots \sigma_{n-1}\cdot w=
\tau^{\alpha+1}\sigma_1\sigma_2\dots \sigma_{n-1}\cdot\sigma_{j-1}^{-1}\dots\sigma_1^{-1}w
$$
and $\sigma_{j-1}^{-1}\dots\sigma_1^{-1} w\in\tilde{W}$. For the third term:
$$\sigma_j^{-1}\dots\sigma_1^{-1}\cdot\tau^{\alpha}\sigma_1\tau\sigma_2\dots \sigma_{n-1}\cdot w=
\sigma_j^{-1}\dots\sigma_1^{-1}\cdot\tau^{\alpha}\sigma_1\sigma_2\dots \sigma_{n-1}\cdot\tau w
$$
and $\tau w\in\tilde{W}$.

\vskip .2cm
(vi) Let now $m=\infty$. We multiply the element (\ref{normalform2}) from the left by $\tau^{-1}$. If $j=0$, there is nothing to prove.
When we multiply the element (\ref{normalform2}), $j>0$, from the left by $\tau^{-1}$, we follow the same steps as in (v), using the identity 
\begin{equation}\label{eq-lem3b}\tau^{-1}\sigma_1^{-1}\tau^{\alpha}\sigma_1=(q-q^{-1})\bigl(\tau^{\alpha}\sigma_1\tau^{-1}-\sigma_1\tau^{\alpha-1}\bigl)+\sigma_1^{-1}\tau^{\alpha}\sigma_1\tau^{-1}\ ,\end{equation}
which is obtained from (\ref{eq-lem3}) by multiplying by $\tau^{-1}$ from both sides.\hfill$\square$

\begin{lemm}
{\hspace{-.2cm}.\hspace{.2cm}}\label{lem3}For any integer $\alpha$, 
we have:
\begin{equation}\label{eq-lem3}\tau\sigma_1^{-1}\tau^{\alpha}\sigma_1=(q-q^{-1})\bigl(\tau\sigma_1\tau^{\alpha}-\tau^{\alpha+1}\sigma_1\bigl)+\sigma_1^{-1}\tau^{\alpha}\sigma_1\tau\ .\end{equation}
\end{lemm}

\emph{Proof of the Lemma.} Multiplying the equality $\sigma_1\tau\sigma_1\tau^{\alpha} =\tau^{\alpha}\sigma_1\tau\sigma_1$ by  $\sigma_1^{-1}$ from 
both sides, we get
\begin{equation}\label{nach2}\tau\sigma_1\tau^{\alpha}\sigma_1^{-1}=\sigma_1^{-1}\tau^{\alpha}\sigma_1\tau\ .\end{equation}
Expanding $\tau\sigma_1^{-1}\tau^{\alpha}\sigma_1=\tau \bigr(\sigma_1-(q-q^{-1})\bigl)\tau^{\alpha} \bigr(\sigma_1^{-1}+(q-q^{-1})\bigl)$ and using 
(\ref{nach2}) we obtain (\ref{eq-lem3}). \hfill$\square$

\subsection{{\hspace{-0.50cm}.\hspace{0.50cm}}Construction of an $H(m,1,n)$-module}

Consider an arbitrary left $H(m,1,n-1)$-module $M_{n-1}$. We denote its elements by letters $u,v$ etc. Let $E_{\infty}:=\mathcal{A}_{\infty}[z,z^{-1}]$ and, for finite $m$,  
$E_m:=\mathcal{A}_m[z]/\langle \chi\rangle$, where $\langle \chi\rangle$ is the ideal generated by the 
characteristic polynomial $\chi(z):=(z-v_1)\dots(z-v_m)$
of $\tau$  (note that $z$ is invertible in $E_m$). 
Let $V$ be a free $\mathcal{A}_m$-module with the basis $w_j$, $j=0,\dots,n-1$. Let $M_n:=V\otimes E_m\otimes M_{n-1}$ (the tensor product is over $\mathcal{A}_m$). The elements 
$w_j\otimes \phi\otimes u$ we denote by ${\cal{V}}_{j,\phi,u}$.
For brevity we write
\begin{equation}\label{def-betaj}\beta_j:=\sigma_{j-1}^{-1}\dots\sigma_{1}^{-1}\ ;\end{equation}
the result of the action of an element $\psi\in H(m,1,n-1)$ on an element $u\in M_{n-1}$ we denote simply by $\psi u$ (without any symbol for the action). 

\vskip .2cm
Define operators $F_{\sigma_i}$, $i=1,\dots,n-1$, and $F_{\tau}$ on the space $M_n$ by (below the last index of ${\cal{V}}$ carries information about the 
action of the elements of the algebra $H(m,1,n-1)$ on the module $M_{n-1}$; for finite $m$, $\phi$ is a polynomial in $z$, defined modulo 
$\chi$, and therefore the element $\phi(\tau)\in H(m,1,n-1)$ which appears in the last index of ${\cal{V}}$ is well defined):
\begin{equation}\label{actionsigma}
F_{\sigma_i}\ : {\cal{V}}_{j,\phi,u}\ \mapsto\ \left\{\begin{array}{ll}
{\cal{V}}_{j,\phi,\sigma_{i-1}u}\ ,& j<i-1\ ,\\[1em]
(q-q^{-1})\, {\cal{V}}_{i-1,\phi,u}+{\cal{V}}_{i,\phi,u}\ ,& j=i-1\ ,\\[1em]
{\cal{V}}_{i-1,\phi,u}\ ,& j=i\ ,\\[1em]
{\cal{V}}_{j,\phi,\sigma_i u}\ ,& j>i\ ,
\end{array}\right.
\end{equation}
and
\begin{equation}\label{actiontau}
F_{\tau}\ :  \left\{\begin{array}{ll}
{\cal{V}}_{0,\phi,u}\ \mapsto\ {\cal{V}}_{0,z\phi,u}\ ,& \\[1em]
{\cal{V}}_{j,\phi,u}\ \mapsto\ (q-q^{-1})\, {\cal{V}}_{0,z,\beta_j\phi(\tau) u}-
(q-q^{-1})\, {\cal{V}}_{0,z\phi,\beta_j u}
+{\cal{V}}_{j,\phi,\tau u}\ ,& j>0\ .
\end{array}\right.
\end{equation}
It follows that $F_{\tau}$ is invertible. 

Let, as above, $\tilde{W}$ denote the subalgebra of the algebra $H(m,1,n)$ generated by the elements $\tau,\tau^{-1}$ and $\sigma_1,\dots,\sigma_{n-2}$.
Take $\tilde{W}$ for the $H(m,1,n-1)$-module $M_{n-1}$ (for the moment we know only that the algebra $\tilde{W}$ is a quotient of $H(m,1,n-1)$; we define the action of $H(m,1,n-1)$
by the left multiplication on its quotient). With the calculation made in the proof of Proposition \ref{normalform-H}, 
one checks that the formulas (\ref{actionsigma}) and (\ref{actiontau}) are valid if one makes the following substitution:
$${\cal{V}}_{j,\phi,u}\ \leadsto\  \sigma_j^{-1}\sigma_{j-1}^{-1}\dots \sigma_1^{-1}\phi(\tau)\sigma_1\sigma_2\dots \sigma_{n-1}u\,,$$ 
and considers $F_{\tau},F_{\sigma_i}$ as operators of the
left multiplication by the corresponding generators. 
More is true. The formulas (\ref{actionsigma}) and (\ref{actiontau}) have the following universal property.

\begin{prop}
{\hspace{-.2cm}.\hspace{.2cm}}\label{normalform-H'} 
The map $\sigma_i\mapsto F_{\sigma_i}$, $i=1,\dots,n-1$, and $\tau\mapsto F_{\tau}$ (and $\tau^{-1}\mapsto F_{\tau}^{-1}$ for $m=\infty$) equips $M_n$ with a structure of an $H(m,1,n)$-module.
\end{prop}

\emph{Proof.} A straightforward although lengthy calculation. Given a defining relation from the list (\ref{rel1a})--(\ref{rel1''b}) 
we verify it on each vector ${\cal{V}}_{j,\phi,u}$. Below we mention different placements of the index $j$ in a verification of a given relation.

\vskip .2cm
(i) For the relation $\sigma_i\sigma_k=\sigma_k\sigma_i$ with $i<k-1$ one considers separately the following positions of the index $j$:
$$j<i-1\ ,\  j=i-1\ ,\  j=i\ ,\  i<j<k-1\ ,\  j=k-1\ ,\  j=k\ \  \textrm{and}\ \  j>k\ .$$

\vskip .2cm
(ii) For the Artin relation $\sigma_i\sigma_{i+1}\sigma_i=\sigma_{i+1}\sigma_i\sigma_{i+1}$ one considers separately the following positions of the index $j$:
$$j<i-1\ ,\  j=i-1\ , \ j=i\ , \ j=i+1\ \  \text{and} \ \ j>i+1\ .$$

\vskip .2cm
(iii) For the relation $\sigma_i^2-(q-q^{-1}) \sigma_i -1=0$ one considers separately the following positions of the index $j$:
$$j<i-1\ , \ j=i-1\ , \ j=i\ \  \text{and}\ \  j>i\ .$$

\vskip .2cm
(iv) For the relation $\tau\sigma_i=\sigma_i\tau$ with $i>1$  one considers separately the following positions of the index $j$:
$$j=0\ , \ j<i-1\ , \ j=i-1\ , \ j=i\ \  \text{and}\ \  j>i\ .$$

\vskip .2cm
(v) For the relation $\,\tau\sigma_1\tau\sigma_1=\sigma_1\tau\sigma_1\tau\,$ it is enough to consider separately the following positions of the index $j$: 
$$j=0\ , \ j=1\ \  \text{and}\ \ j>1\ .$$
The following observation simplifies a verification here:
\begin{equation}F_{\tau}F_{\sigma_1}F_{\tau}:  \left\{\!\!\!\begin{array}{lcl}
{\cal{V}}_{0,\phi,u}&\!\!\!\mapsto\!\!\!&(q-q^{-1})\, {\cal{V}}_{0,z,\tau\phi(\tau)u}+{\cal{V}}_{1,z\phi,\tau u}\ , \\[1.5em]
{\cal{V}}_{1,\phi,u}&\!\!\!\mapsto\!\!\!&(q-q^{-1})\, {\cal{V}}_{1,z,\tau\phi(\tau) u}-(q-q^{-1})\, {\cal{V}}_{1,z\phi,\tau u}
+{\cal{V}}_{0,z\phi,\tau u}\ ,  \\[1.5em]
{\cal{V}}_{j,\phi,u}&\!\!\!\mapsto\!\!\!&(q-q^{-1})^2\, {\cal{V}}_{0,z,\tau [ \beta_j,\phi(\tau)] u}+
(q-q^{-1})\, {\cal{V}}_{1,z,\tau \beta_j\phi(\tau) u}-(q-q^{-1})\, {\cal{V}}_{1,z\phi,\tau \beta_j u} \\[1em]
&&+(q-q^{-1})\, {\cal{V}}_{0,z,\beta_j \phi(\tau)\sigma_1\tau u}-(q-q^{-1})\, {\cal{V}}_{0,z\phi,\beta_j \sigma_1\tau u}
+{\cal{V}}_{j,\phi,\tau\sigma_1\tau u}\ , \ \ j>1\ .
\end{array}\right.
\end{equation}
Here $ [ \beta_j,\phi(\tau)]$ is the commutator of $\beta_j$ and $\phi(\tau)$; in the verification of the relation  
$\tau\sigma_1\tau\sigma_1=\sigma_1\tau\sigma_1\tau$ on ${\cal{V}}_{j,\phi,u}$ with $j>1$ we use the formula 
(\ref{eq-lem3}) in the form
$$\tau\sigma_1^{-1}\phi(\tau)\sigma_1=(q-q^{-1})\bigl(\tau\sigma_1\phi(\tau)-\tau\phi(\tau)\sigma_1\bigl)+\sigma_1^{-1}\phi(\tau)\sigma_1\tau\ .$$

\vskip .2cm
(vi) Let $m<\infty$. For the relation $(\tau-v_1)\dots (\tau-v_m)=0$ one considers separately the following positions of the index $j$: 
$$j=0\ \  \text{and}\ \  j>0\ .$$ 
A verification of this relation for $j>0$ is the only place in the proof which maybe requires a comment. For $j>0$ one proves by induction the 
following formula:
\begin{equation}\label{pui-tau}F_{\tau}^l\ :  {\cal{V}}_{j,\phi,u}  \mapsto (q-q^{-1}) \sum\limits_{i=1}^l {\cal{V}}_{0,z^i,\beta_j\phi(\tau)\tau^{l-i}u}-(q-q^{-1}) 
\sum\limits_{i=1}^l {\cal{V}}_{0,z^i\phi,\beta_j\tau^{l-i}u}+{\cal{V}}_{j,\phi,\tau^l u}\ ,\ \  j>0\ .\end{equation}
The first sum in (\ref{pui-tau}) can be seen as the image of an element 
\begin{equation}\label{imafs}(z\otimes \phi )\cdot\frac{1\otimes z^l-z^l\otimes 1}{1\otimes z -z\otimes 1}\end{equation} 
from the space
$E_m\otimes E_m$ to $M_n$ with respect to the map $\kappa_u$, defined for each $u\in M_{n-1}$ by 
$$\kappa_u\, :\ f\otimes g\mapsto  {\cal{V}}_{0,f,\beta_j g(\tau )u}\ .$$
In (\ref{imafs}), the fraction $\frac{1\otimes z^l-z^l\otimes 1}{1\otimes z -z\otimes 1}$ is understood as the image of a polynomial which is the result of 
the division of the numerator by the denominator (as polynomials of two unrestricted variables) in the space $E_m\otimes E_m$. Similarly, the second sum in  
(\ref{pui-tau}) can be understood as the image of $(z\phi\otimes 1)\cdot\frac{1\otimes z^l-z^l\otimes 1}{1\otimes z -z\otimes 1}$ 
with respect to the same 
map $\kappa_u$. Thus the  first sum minus the second sum (the combination which appears in  (\ref{pui-tau})) is the image of 
\begin{equation}\label{imafs2}(1\otimes \phi- \phi\otimes 1)\cdot (z\otimes 1)\cdot\frac{1\otimes z^l-z^l\otimes 1}{1\otimes z -z\otimes 1}\ .\end{equation} 
The element (\ref{imafs2}) already as a polynomial (and therefore as an element in $E_m\otimes E_m$) can be written in the form
\begin{equation}\label{imafs2'}\frac{1\otimes \phi- \phi\otimes 1}{1\otimes z -z\otimes 1}\cdot (z\otimes 1)\cdot (1\otimes z^l-z^l\otimes 1)\ ,\end{equation}
where the fraction $\frac{1\otimes \phi- \phi\otimes 1}{1\otimes z -z\otimes 1}$ is understood again as the image of a polynomial which is the result of 
the division of the numerator by the denominator (as polynomials of two unrestricted variables) in the space $E_m\otimes E_m$. 

\vskip .2cm
Writing now $\chi (z)=\sum_{l=0}^m c_m z^m$ one verifies the relation $\chi (F_{\tau})=0$ immediately with the help of (\ref{imafs2'})
(recall that $E_m=\mathcal{A}_m[z]/\langle\chi\rangle$). \hfill$\square$

\begin{rema} {\hspace{-0.2cm}.\hspace{0.2cm}} {\rm The action of $F_{\tau}$ on the vectors of the form ${\cal{V}}_{j,1,u}$ with $j>0$ is simply the action of $\tau$ on $M_{n-1}$, that is,
\begin{equation}\label{actiontau2}
F_{\tau}\ :\ {\cal{V}}_{j,1,u}\ \mapsto\ {\cal{V}}_{j,1,\tau u}\ \ \ {\text{for}}\ \  j>0\ .
\end{equation}
\hfill$\triangle$ 
}\end{rema}

\begin{rema}{\hspace{-0.2cm}.\hspace{0.2cm}}  {\rm The operators $F_{\sigma_i}$ and $F_{\tau}$ defined in (\ref{actionsigma}) and (\ref{actiontau}) can be represented by 
$n\times n$ matrices (with indices related to the basis of $V$) whose elements are operators acting in the space $E_m\otimes M_{n-1}$. 
By $\hat{z}$ we denote the operator of the multiplication by $z$ in the space $E_m$. To fit formulas in the line,
we denote the operator $\text{Id}_{E_m}\otimes \sigma_i$ simply by $\sigma_i$, the operator $\text{Id}_{E_m}\otimes \tau$ simply by $\tau$  and the operator  
$\hat{z}\otimes \text{Id}_{M_{n-1}}$ simply by $\hat{z}$.

\vskip .2cm
The operator $F_{\sigma_i}$ reads (recall that the elements of the basis of $V$ are labelled by numbers from 0 to $n-1$)
\begin{equation}\label{mfoofs}
F_{\sigma_i}=\left(\begin{array}{ccccccc}\sigma_{i-1}&&&&&& \\ &\ddots&&&&& \\ &&\sigma_{i-1}&&&&\\ &&&q-q^{-1}&1&&\\[.5em] &&&1&0&&\\ 
&&&&&\sigma_i&\\ &&&&&&\ddots
\end{array}\right)\ ;\end{equation}
\vskip .2cm\noindent
here the $2\times 2$ block with ones (that is, the identity operators) outside the main diagonal is in the 
$(i-1)^{\text\footnotesize{st}}$ and $i^{\text\footnotesize{th}}$ lines and columns. 

\vskip .2cm
The operator $F_{\tau}$ reads  
\begin{equation}\label{mfooft}
F_{\tau}=\left(\begin{array}{cccccc}\hat{z}&\text{\small{$(q-q^{-1})$}}\hat{z}\text{\footnotesize{$(\mu -1)$}}&
\text{\footnotesize{$(q-q^{-1})$}}\hat{z}\sigma_1^{-1}\text{\footnotesize{$(\mu -1)$}}&
\text{\footnotesize{$(q-q^{-1})$}}\hat{z}\sigma_2^{-1}\sigma_1^{-1}\text{\footnotesize{$(\mu -1)$}}&&\dots \\[1.5em] 
&\tau&&&&\dots \\[1.5em] &&\tau&&&\dots\\[1.5em] &&&\tau&&\dots\\[1.5em] 
\vdots&\vdots&\vdots&\vdots&&\ddots
\end{array}\right)\ ;
\end{equation}
here only the first line and the main diagonal have non-zero entries. The operator $\mu$ on the space $E_m\otimes M_{n-1}$ is defined as follows:
$$\mu (\phi\otimes u):=1\otimes \phi(\tau)u\ ,$$
where $\phi$ is a polynomial in $z$. The operator $\mu$ has the following properties:
$$\mu\hat{z}=\tau\mu\ ,\  \mu\tau =\tau\mu\ ,\ \mu^2=\mu\ .$$
\hfill$\triangle$ 
}\end{rema}

\section{\hspace{-0.55cm}.\hspace{0.55cm}Flatness of deformation. Normal form for elements of $H(m,1,n)$}\label{sec-nf2}

We are now ready to prove that the deformation $H(m,1,n)$ of the group ring $\mathbb{C}G(m,1,n)$ is flat and to give the normal form for 
elements of the algebra $H(m,1,n)$.

Some results of this section were used in \cite{OPdA2} and \cite{OPdA}.

\subsection{{\hspace{-0.50cm}.\hspace{0.50cm}}Flatness of deformation}

As above, $\tilde{W}$ denotes the subalgebra of the algebra $H(m,1,n)$ generated by the elements $\tau,\tau^{-1}$ 
and $\sigma_1,\dots,\sigma_{n-2}$. 

\begin{prop}
{\hspace{-.2cm}.\hspace{.2cm}}\label{dimension} 
{\rm (i)} The subalgebra $\tilde{W}$ is isomorphic to $H(m,1,n-1)$.

\vskip .2cm\noindent
{\rm (ii)} The algebra $H(m,1,n)$ is a free $\mathcal{A}_m$-module. Let $\B^{(n-1)}$ 
be a basis of $H(m,1,n-1)\cong\tilde{W}$. Then the following set of elements 
is a basis of $H(m,1,n)$:
\begin{equation} \sigma_j^{-1}\sigma_{j-1}^{-1}\dots \sigma_1^{-1}\tau^{\alpha}\sigma_1\sigma_2\dots \sigma_{n-1}w\ , \ j\in\{0,1,\dots ,n-1\}\,, \ \alpha\in 
\mathfrak{E}_m\,,\ w\in \B^{(n-1)} 
\ . \end{equation}
Therefore, the algebra $H(m,1,n)$ is
a flat deformation of the group ring of $G(m,1,n)$. 
\end{prop}

\emph{Proof.} Assume, by induction, that $H(m,1,n-1)$ is a free $\mathcal{A}_m$-module and a flat deformation of the group ring of $G(m,1,n-1)$ (the induction basis is trivial since $H(m,1,0)=\mathcal{A}_m$). Let $\mathfrak{p}\ :\ H(m,1,n-1)\to\tilde{W}$ denote the natural homomorphism. The Proposition \ref{normalform-H} implies that the elements
\begin{equation}\label{normalform2b}
\sigma_j^{-1}\sigma_{j-1}^{-1}\dots \sigma_1^{-1}\tau^{\alpha}\sigma_1\sigma_2\dots \sigma_{n-1}\mathfrak{p}(w)\ ,
\end{equation}
where $j\in\{0,\dots,n-1\}$, $\alpha\in\mathfrak{E}_m$ and $w\in \B^{(n-1)}
$,
span $H(m,1,n)$ as a vector space. We shall show that these elements are linearly independent.

Let $M_{n-1}$ be the left regular module for $H(m,1,n-1)$; that is, the module space is the algebra itself and the elements of the algebra act by left multiplication. 

By the Proposition \ref{normalform-H'}, the space $M_n=V\otimes E_m\otimes M_{n-1}$ is equipped 
with an $H(m,1,n)$-module structure. Denote by $F_a$ the operator corresponding to the element $a$ for any $a\in H(m,1,n)$. Using formulas (\ref{actionsigma}) and  (\ref{actiontau}) for the action of $H(m,1,n)$ on $M_n$, and also formula (\ref{actiontau2}), we obtain
\begin{equation}F_{\sigma_j^{-1}\dots\sigma_1^{-1}\phi(\tau)\sigma_1\dots\sigma_{n-1}\mathfrak{p}(w)}\ :\ \ {\cal{V}}_{n-1,1,1}\ \mapsto\ {\cal{V}}_{j,\phi,w}\ .\end{equation} 
Thus, the operators 
$F_{\sigma_j^{-1}\dots\sigma_1^{-1}\phi(\tau)\sigma_1\dots\sigma_{n-1}\mathfrak{p}(w)}$, $j\in\{0,\dots,n-1\}$, $\alpha\in\mathfrak{E}_m$ and $w\in
\B^{(n-1)}
$, are linearly independent, which implies the linear independence of the set (\ref{normalform2b}). The non-triviality of the kernel of $\mathfrak{p}$ would contradict to 
the linear independence of the elements (\ref{normalform2b}). Therefore, the subalgebra $\tilde{W}$ is isomorphic to $H(m,1,n-1)$ and the flatness of the deformation from the group ring of $G(m,1,n)$ to $H(m,1,n)$ follows. \hfill$\square$

\subsection{{\hspace{-0.50cm}.\hspace{0.50cm}}Basis $\mathcal{B}$ of $H(m,1,n)$}

We construct recursively, in the same way as we did for $G(m,1,n)$, a global normal form for elements of $H(m,1,n)$ using now Proposition \ref{dimension}, statement (ii). Let $\B_k$ 
be the set of elements 
$\{\sigma_j^{-1}\sigma_{j-1}^{-1}\dots \sigma_1^{-1}\tau^{\alpha}\sigma_1\sigma_2\dots \sigma_{k-1}\ \vert\ j\in\{0,\dots,k-1\}\,,\ \alpha\in\mathfrak{E}_m\}$.

\begin{coro}
{\hspace{-.2cm}.\hspace{.2cm}}\label{normalform-fin} 
Any element $x\in H(m,1,n)$ can be written uniquely as a linear combination of elements  
\begin{equation}\label{nfch}x=u_nu_{n-1}\dots u_1\ ,\end{equation} 
where $u_k\in \B_k$ for $k=1,\dots,n$. 
\end{coro}

In other words, the products $u_nu_{n-1}\dots u_1$, where $u_k$ ranges over $\B_k$ 
for $k=1,\dots,n$, form a basis of the $\mathcal{A}_m$-module $H(m,1,n)$. 
We denote such basis, for brevity, $\mathcal{B}:=\B_{n}\dots \B_{1}$. 

\begin{rema}{\hspace{-0.2cm}.\hspace{0.2cm}} \label{rem-H}
{\rm Let $\mathcal{A}_mH_n:=\mathcal{A}_m\otimes H_n$, where the tensor product is over $\mathbb{C}[q,q^{-1}]$ (recall that $H_n$ is the Hecke algebra of type A, see Section \ref{sec-def}).
Define the homomorphism $\varsigma\ :\ \mathcal{A}_mH_n\rightarrow H(m,1,n)$ 
by sending the generator $\sigma_i$ of $\mathcal{A}_mH_n$ to the generator $\sigma_i$ of $H(m,1,n)$ for $i=1,\dots ,n-1$. 
Let $\B_{k,0}$ be the subset of $\B_k$ 
of all elements with $\alpha=0$. The image of the standard basis of $\mathcal{A}_mH_n$ consists of products (\ref{nfch}) where 
$u_k\in \B_{k,0}$ for $k=1,\dots,n$. Thus, by the linear independence of the words (\ref{nfch}), $\varsigma$ is injective.

There is another way, without the use of the Corollary \ref{normalform-fin}, to check that $\varsigma$ is an embedding. If $m$ is finite, 
fix a number $e$, $1\leq e\leq m$. The map, which sends the generator $\sigma_i$ of $H(m,1,n)$ to the generator $\sigma_i$ of $\mathcal{A}_mH_n$ for $i=1,\dots ,n-1$
and the generator $\tau$ of $H(m,1,n)$ to the element $v_e$, clearly extends to a homomorphism $\pi
: H(m,1,n)\rightarrow \mathcal{A}_mH_n$ (if $m=\infty$, the element $v_e$ is an arbitrary unit in $\mathcal{A}_{\infty}$). One has 
$\pi\circ\varsigma =\text{Id}_{\mathcal{A}_mH_n}$ so $\pi$ is a left inverse to $\varsigma$
; in particular, $\varsigma$ is an embedding. Likewise, we have, 
on the level of braid groups, an embedding $B_n\to\alpha B_n$ defined 
by a map, tautological on the generators of $B_n$.
\hfill$\triangle$ 
}\end{rema}

\subsection{{\hspace{-0.50cm}.\hspace{0.50cm}}Induced representations}

Let  ${\mathfrak{B}}$ be an associative subalgebra of an associative algebra ${\mathfrak{A}}$. Let $W$ be a left ${\mathfrak{B}}$-module. The vector space $ {\mathfrak{A}}\otimes_ {\mathfrak{B}} W$ carries a natural ${\mathfrak{A}}$-module structure 
defined by ${\mathfrak{a}}.({\mathfrak{a}}'\otimes w):={\mathfrak{a}}{\mathfrak{a}}'\otimes w$. This is the induced ${\mathfrak{A}}$-module.

\vskip .2cm 
Let $M_{n-1}$ be a left $H(m,1,n-1)$-module and  $M_n:=V\otimes E_m\otimes M_{n-1}$ 
the $H(m,1,n)$-module, constructed in the preceding Subsection. 
It follows from the Proposition \ref{dimension} that $M_n$ is the induced $H(m,1,n)$-module with respect to the subalgebra $H(m,1,n-1)$ and the module $M_{n-1}$ over it.  
The formulas (\ref{actionsigma}) and (\ref{actiontau}) give an explicit realization of the induced module $M_n$.
In particular, the construction (\ref{actionsigma})--(\ref{actiontau}) possesses the following functoriality properties. For a 
free $\mathcal{A}_m$-module $M$, let $\Upsilon (M):=V\otimes E_m\otimes M$. For two free $\mathcal{A}_m$-modules $M,M'$ and a map $\alpha :M\rightarrow M'$, let 
$\Upsilon (\alpha):=\text{Id}_V\otimes \text{Id}_{E_m}\otimes \alpha$. Then $\Upsilon$ is a functor from the category of $H(m,1,n-1)$-modules 
to the category of $H(m,1,n)$-modules,
$$\Upsilon \,:\, H(m,1,n-1)\text{-mod}\ \rightarrow\ H(m,1,n)\text{-mod}\ .$$ 
The formulas (\ref{actionsigma})--(\ref{actiontau}) do not contain denominators, so specializations (setting $q$, $v_1,\dots,v_m$ to 
arbitrary complex numbers different from 0 in the formulas (\ref{actionsigma})--(\ref{actiontau})) do not cause any difficulties. This concerns,  in particular, the specialization 
(\ref{specialization}) to the group ring.

\begin{rema}{\hspace{-0.2cm}.\hspace{0.2cm}}{\rm
The space $E_1$ is one-dimensional and we identify it with the ring $\mathcal{A}_1$.
We have here the operators $F_{\sigma_i}$ acting on the space $V\otimes M_{n-1}$:
\begin{equation}\label{actionsigma2}
F_{\sigma_i}\ : {\cal{V}}_{j,u}\ \mapsto\ \left\{\begin{array}{ll}
{\cal{V}}_{j,\sigma_{i-1}u}\ ,& j<i-1\ ,\\[1em]
(q-q^{-1})\, {\cal{V}}_{i-1,u}+{\cal{V}}_{i,u}\ ,& j=i-1\ ,\\[1em]
{\cal{V}}_{i-1,u}\ ,& j=i\ ,\\[1em]
{\cal{V}}_{j,\sigma_i u}\ ,& j>i\ .
\end{array}\right.\end{equation}
The formula (\ref{actionsigma2}) equips the space $M_n:=V\otimes M_{n-1}$ with an $\mathcal{A}_1 H_n$-module structure. Since the formula (\ref{actionsigma2})
does not involve the parameter $v_1$, we have, in fact, an $H_n$-module structure on $V\otimes M_{n-1}$ where $M_{n-1}$ is a left $H_{n-1}$-module and 
$V$ is understood as a free $\mathbb{C}[q,q^{-1}]$-module. So, for $m=1$, $\Upsilon$ is a functor from the category of $H_{n-1}$-modules to the category of $
H_{n}$-modules,
$$\Upsilon \,:\, H_{n-1}\text{-mod}\ \rightarrow\ H_{n}\text{-mod}\ .$$  
The label $\phi$ of vectors ${\cal{V}}_{j,\phi,u}$ is not touched by the action of the operators $F_{\sigma_i}$ given by the formula (\ref{actionsigma}).
Thus the restriction of the $H(m,1,n)$-module $M_n$ to the subalgebra generated by $\sigma_1,\dots,\sigma_{n-1}$ (isomorphic to $\mathcal{A}_m H_n$, see Remark \ref{rem-H}, is the direct sum of $m$ copies of the $\mathcal{A}_m H_n$-module given by the formula (\ref{actionsigma2}).
\hfill$\triangle$ 
}\end{rema}

\begin{rema}{\rm
{\bf (Burau module).} Take for $M_{n-1}$ the one-dimensional $H_{n-1}$-module, on which the generators $\sigma_i$ act by multiplication by 
$q$. Then the resulting module $M_n$ is the Burau module for the Hecke algebra $H_n$.

Taking now for $M_{n-1}$ the one-dimensional module of the algebra $H(m,1,n-1)$ on which the generators $\sigma_i$ 
act by multiplication by $q$ and the generator $\tau$ acts by multiplication by $v_e$ for some $e$, $1\leq e\leq m$,
(if $m=\infty$, $v_e$ is an arbitrary unit in $\mathcal{A}_{\infty}$), the resulting  $H(m,1,n)$-module $M_n\cong V\otimes E_m$ is 
an analogue of the Burau module. The action of the generators $\sigma_i$ is given by the usual Burau matrices (these are 
the matrices (\ref{mfoofs}) in which every $\sigma$ is replaced by $q$; these matrices act trivially in the space $E_m$) 
while the matrix of the operator $\tau$ 
is given by
\begin{equation}\label{mfooftm}  
F_{\tau}=\left(\begin{array}{cccccc}\hat{z}&\text{\small{$(q-q^{-1})$}}\hat{z}\text{\footnotesize{$(\mu_e -1)$}}&
\text{\footnotesize{$(q-q^{-1})$}}q^{-1}\hat{z}\text{\footnotesize{$(\mu_e -1)$}}&
\text{\footnotesize{$(q-q^{-1})$}}q^{-2}\hat{z}\text{\footnotesize{$(\mu_e -1)$}}&&\dots \\[1.5em] 
&v_e&&&&\dots \\[1.5em] &&v_e&&&\dots\\[1.5em] &&&v_e&&\dots\\[1.5em]
\vdots&\vdots&\vdots&\vdots&&\ddots
\end{array}\right)\ ;
\end{equation}
here $\mu_e$ is defined by $\mu_e(\phi):=\phi (v_e)$, where $\phi\in E_m$.
\hfill$\triangle$ 
}\end{rema}

\subsection{{\hspace{-0.50cm}.\hspace{0.50cm}}Other bases of $H(m,1,n)$}\label{sub-otherbases}

We mention several other normal forms for elements of $H(m,1,n)$ with respect to $\tilde{W}\cong H(m,1,n-1)$ similar to the form from the Proposition \ref{normalform-H} (and thus several other inductive bases of $H(m,1,n)$ similar to the basis $\B$ in the Corollary \ref{normalform-fin}). 

\vskip .2cm
Let, as above, $\tilde{W}$ be the subalgebra generated by $\tau,\tau^{-1},\sigma_1,\dots,\sigma_{n-2}$. An induction on $j$ establishes the following fact.
\begin{lemm}{\hspace{-.2cm}.\hspace{.2cm}}\label{auxile} 
Let $\epsilon =\pm 1$. An element of the form 
$$\sigma_j^{\epsilon}\sigma_{j-1}^{\epsilon}\dots \sigma_1^{\epsilon}\tau^{\alpha}\sigma_1\sigma_2\dots \sigma_{n-1} w\ \ \text{with $j\in\{0,\dots,n-1\}$, $\alpha\in
\E_m$ and $w\in \tilde{W}$}$$
equals to $\sigma_j^{-\epsilon}\sigma_{j-1}^{-\epsilon}\dots \sigma_1^{-\epsilon}\tau^{\alpha}\sigma_1\sigma_2\dots \sigma_{n-1} w$ plus a linear combination of elements 
$$\sigma_k^{-\epsilon}\sigma_{k-1}^{-\epsilon}\dots \sigma_1^{-\epsilon}\tau^{\alpha}\sigma_1\sigma_2\dots \sigma_{n-1} w_k\ \ \text{with $k<j$ and $w_k\in \tilde{W}$}\ .$$ 
\end{lemm}

Therefore, any $x\in H(m,1,n)$ can be written as a linear combination of elements of the set 
\begin{equation}\label{normalform2'}\left\{\begin{array}{l}
\sigma_{k+1}\sigma_{k+2}\dots \sigma_{n-1}w\ ,\ \textrm{where $k\in\{0,\dots,n-1\}$ and $w\in \tilde{W}$\ ,}\\[.5cm]
\sigma_j\sigma_{j-1}\dots \sigma_1\tau^{\alpha}\sigma_1\dots \sigma_{n-1}w\ ,\ \textrm{where $j\in\{0,\dots,n-1\}$, $\alpha\in\mathfrak{E}_m\backslash\{0\}
$ and $w\in \tilde{W}$\ .}\end{array}\right. 
\end{equation}

Let $\iota$ be the involution of the ring $H(m,1,n)$, defined by
\begin{equation}\label{iota}
\iota(x)=x^{-1},\ x\in\{\tau,\sigma_1,\dots,\sigma_{n-1}\},\ \quad\iota(q)=q^{-1}\ \quad\text{and, for $m<\infty$,}\ \ \iota(v_a)=v_a^{-1},\ a=1,\dots,m.
\end{equation}
Applying the involution $\iota$ to the normal forms (\ref{normalform2}) and  (\ref{normalform2'}), we find that any $x\in H(m,1,n)$ can be written as a linear combination of elements of the set 
\begin{equation}\label{normalform2''}
\sigma_j\sigma_{j-1}\dots \sigma_1\tau^{\alpha}\sigma_1^{-1}\dots \sigma_{n-1}^{-1}w\ ,\ \textrm{where $j\in\{0,\dots,n-1\}$, $-\alpha\in\mathfrak{E}_m$ and $w\in \tilde{W}$\ ,}
\end{equation}
or of the set
\begin{equation}\label{normalform2'''}
\!\left\{\!\!\begin{array}{l}
\!\sigma^{-1}_{k+1}\sigma^{-1}_{k+2}\dots \sigma^{-1}_{n-1}w\ ,\ \textrm{where $k\in\{0,\dots,n-1\}$ and $w\in \tilde{W}$\ ,}\\[1em]
\!\sigma_j^{-1}\sigma_{j-1}^{-1}\dots \sigma_1^{-1}\tau^{\alpha}\sigma_1^{-1}\dots \sigma_{n-1}^{-1}w\ ,\ 
\textrm{where $j\in\{0,\dots,n-1\}$, $-\alpha\in\mathfrak{E}_m\backslash\{0\}$ and $w\in \tilde{W}$\ .}
\end{array}\right.
\end{equation}

Four other normal forms for elements of $H(m,1,n)$ are obtained by the application of the anti-involution $\varpi$, $\varpi(xy)=\varpi(y)\varpi(x)$, of $H(m,1,n)$, identical on the generators  $\tau,\tau^{-1},\sigma_1,\dots,\sigma_{n-1}$, to the normal forms (\ref{normalform2}), (\ref{normalform2'}), (\ref{normalform2''}) and (\ref{normalform2'''}).

\section{\hspace{-0.55cm}.\hspace{0.55cm} Relative traces for the chain of algebras $H(m,1,n)$ 
}\label{ConditExpec}

We examine and prove the existence and uniqueness of linear maps $\Tr_k$ from the algebra $H(m,1,k)$ to the algebra $H(m,1,k-1)$,
$k=1,2,\dots$, 
which satisfiy
\begin{equation}\label{cond-exp-val1}
\Tr_1(1)=1\ \ \ \text{and}\ \ \ \Tr_1(\tau^a)=\mu_a\ \text{where}\ 
\ a\in \E_m\setminus\{0\}\ ,\ \mu_a\in 
\mathcal{A}_m\ ,
\end{equation}
and, for $k\geq2$, $X,Y\in H(m,1,k-1)$ and $Z\in H(m,1,k)$,
\begin{eqnarray}
\label{cond-exp-val2a}
&\Tr_k(XZY)=X\Tr_k(Z)Y\ ,&\\[0.4em]
\label{cond-exp-val2b}
&\Tr_k\left(\sigma_{k-1}^{\varepsilon}X\sigma_{k-1}^{-\varepsilon}\right)=\Tr_{k-1}(X)\ \ \ \text{where $\varepsilon=\pm1$\ ,}&\\[0.4em]
\label{cond-exp-val2c}
&\Tr_{k-1}(\Tr_k(\sigma_{k-1}Z))=\Tr_{k-1}(\Tr_k(Z\sigma_{k-1}))\ 
\ ,&\\[0.4em]
\label{cond-exp-val2d}
&\Tr_k(\sigma_{k-1})=D\ \ \text{where $D\in 
\mathcal{A}_m$.}&
\end{eqnarray}
Equality (\ref{cond-exp-val2b}), for $X=1$, together with the initial condition (42), implies that $\Tr_k(1)=1$.
The elements $D$ and $\mu_a$, $a\in \E_m\setminus\{0\}$, are the parameters of these linear maps. For later convenience, we set $\mu_0:=1$ and define, for finite $m$, elements $\mu_a\in \mathcal{A}_m$ also for $a\geq m$ by $\mu_a:=\Tr_1(\tau^a)$, $a\geq m$.
One can also consider  $D$ and $\mu_a$, $a\in \E_m\setminus\{0\}$, as indeterminates and work
over the ring 
of polynomials in these indeterminates with coefficients in $\mathcal{A}_m$. 

\vskip .2cm
The \emph{relative traces} $\Tr_k$, $k=1,2,\dots$,  
serve for a construction of commutative subalgebras, containing the Hamiltonian of the chain models, of the cyclotomic and affine Hecke algebras.

\begin{rema}{\hspace{-0.2cm}.\hspace{0.2cm}}  \label{rema-mark}
{\rm Assume that linear maps $\Tr_k$, $k=1,2,\dots$, satisfying the conditions (\ref{cond-exp-val1})--(\ref{cond-exp-val2d}) exist.
Define the following linear function on $H(m,1,n)$ with values in $\mathcal{A}_m$:
$$\Tr:=\Tr_1\circ\dots\circ\Tr_{n-1}\circ\Tr_n\ .$$
Condition (\ref{cond-exp-val2a}), together with $\Tr_k(1)=1$ for any $k$, 
ensure that the map $\Tr$ is compatible with the chain, in
$n$, of algebras $H(m,1,n)$,
in the sense that $\Tr_1\circ\dots\circ\Tr_{n-1}\circ\Tr_n (X)=\Tr_1\circ\dots\circ\Tr_{n-1}(X)$ for $X\in H(m,1,n-1)$.
As a consequence of the conditions (\ref{cond-exp-val1})--(\ref{cond-exp-val2d}), it is straightforward 
to check that the linear map $\Tr$ satisfies the following properties:
\begin{eqnarray}
\label{proptrmara}
&\Tr(1)=1\ ,&\\[0.5em]
\label{proptrmarb}
&\Tr(ZZ')=\Tr(Z'Z)\,,\ \ \, Z,Z'\in H(m,1,n)\ ,&\\[0.5em]
\label{proptrmarc}
&\Tr(\sigma_{n-1}X)=D\,\Tr(X) \ ,\ \ X\in H(m,1,n-1)\ ,&
\end{eqnarray}
and, in addition,
\begin{equation}\label{proptrmard}
\Tr(\sigma_{n-1}\dots\sigma_1\tau^a\sigma_1^{-1}\dots\sigma_{n-1}^{-1}X)=\mu_a\Tr(X)\ ,\ \ a\in\mathfrak{E}_m\,,\ \ \, X\in H(m,1,n-1)\,.
\end{equation}
Conditions (\ref{proptrmara})--(\ref{proptrmarc}) assert 
that $\Tr$ is a Markov trace on the algebra $H(m,1,n)$ with Markov parameter $D$. The additional property (\ref{proptrmard}) shows that $\Tr$ coincides with the Markov trace studied in \cite{L} with parameters $D$ and $\mu_a$, $a\in \E_m\setminus\{0\}$. 
By results of \cite{L}, such a Markov trace exists and is uniquely determined by the elements $D$ and $\mu_a$, $a\in \E_m\setminus\{0\}$. The results in this section give in particular another proof of the existence (the unicity is easy to check with the basis $\mathcal{B}$), and moreover show that every such Markov trace can be obtained as a composition of relative traces. 

The central forms $\iota(L^{\gamma}_{\ n})$ which we will discuss in more detail in the next Section correspond to a choice $D=0$ and $\mu_a=\iota^{(0)}(\gamma_{-a})$, $a\in\E_m\setminus\{0\}$, see also Remark \ref{markov}.
\hfill$\triangle$ 
}\end{rema}

Let $k\in\{1,\dots,n\}$. For $j\in\{0,\dots,k-1\}$, $a\in\mathbb{Z}$ and $u\in H(m,1,k-1)$, define the following element of $H(m,1,k)$:
\begin{equation}\label{def-T}
T^{(k)}_{j,a,u}:=\sigma_j^{-1}\sigma_{j-1}^{-1}\dots\sigma_1^{-1}\tau^a\sigma_1\dots\sigma_{k-2}\sigma_{k-1}u\ .\end{equation}
For an arbitrary basis $\B^{(k-1)}$ of $H(m,1,k-1)$, by Proposition \ref{dimension}, the elements 
\begin{equation}\label{base-T}T^{(k)}_{j,a,u}\ \ \ \text{for $j\in\{0,\dots,k-1\}$, $a\in\mathfrak{E}_m$ and $u\in \B^{(k-1)}
$,}\end{equation}
form a basis of $H(m,1,k)$.
The formulas (\ref{actionsigma}) and (\ref{actiontau}) (with $T^{(n)}_{j,a,u}$ instead of $\mathcal{V}_{j,\phi,u}$)
give the left regular action of the generators of $H(m,1,n)$ in the basis (\ref{base-T}).
It is straightforward to see that conditions (\ref{cond-exp-val1})--(\ref{cond-exp-val2d}) imply 
\begin{equation}\label{def-tr}
\Tr_k(T^{(k)}_{j,a,u})=\left\{\begin{array}{l}
D\sigma_j^{-1}\dots\sigma_1^{-1}\tau^a\sigma_1\dots \sigma_{k-2} 
u\ \ \ \ \text{if $j<k-1$,}\\[0.4em]
\mu_a u\ \ \ \ \text{if $j=k-1$.}
\end{array}\right.
\end{equation}
Thus, if linear maps $\Tr_k$, $k=1,2,\dots$, satisfying the conditions (\ref{cond-exp-val1})--(\ref{cond-exp-val2d}) exist, they are uniquely determined by (\ref{def-tr}) and by linearity. The next proposition shows the existence. 
\begin{prop}\hspace{-.2cm}.\hspace{.2cm}\label{prop-tr}
The linear maps ${\rm{Tr}}_k\ :\ H(m,1,k)\to H(m,1,k-1)$, $k=1,2,\dots$, 
defined by (\ref{def-tr}), satisfy the conditions (\ref{cond-exp-val1})--(\ref{cond-exp-val2d}).
\end{prop}
\emph{Proof.} (i) As $T^{(k)}_{k-1,0,1}=1$ for $k\geq 1$, $T^{(1)}_{0,a,1}=\tau^a$ for $a\in \E_m$ and $T^{(k)}_{k-2,0,1}=\sigma_{k-1}$ for $k\geq 2$, we immediately get from  (\ref{def-tr}) that (\ref{cond-exp-val1}) and (\ref{cond-exp-val2d}) are satisfied. 

\vskip .2cm
(ii) We verify (\ref{cond-exp-val2a}).
It is enough to consider the situation $k>1$, $Z=T^{(k)}_{j,a,u}$ for $j\in\{0,\dots,k-1\}$, 
$a\in \E_m$ and $u\in H(m,1,k-1)$, and $X,Y\in\{\tau, \tau^{-1}, \sigma_1,\dots,\sigma_{k-2}\}$. 

Fix $Y\in\{\tau,\tau^{-1},\sigma_1,\dots,\sigma_{k-2}\}$. An easy calculation involving 
formulas (\ref{actionsigma}) (with $\mathcal{V}_{j,\phi,u}$ replaced by $T^{(k)}_{j,a,u}$) shows that,
for $i\in\{1,\dots,k-2\}$,
\begin{equation}\label{pr-tr1}\Tr_k(\sigma_i T^{(k)}_{j,a,u} Y)=\left\{\begin{array}{ll}
D \sigma_j^{-1}\dots \sigma_1^{-1}\tau^a\sigma_1\dots\sigma_{k-2}\sigma_{i-1}\, u\,Y \hspace{1cm} & \text{if $j<i-1$,}\\[0.5em]
D \sigma_{j+1} \sigma_j^{-1}\dots \sigma_1^{-1}\tau^a\sigma_1\dots\sigma_{k-2}\, u\,Y & \text{if $j=i-1$,}\\[0.5em]
D \sigma_{j-1}^{-1}\dots \sigma_1^{-1}\tau^a\sigma_1\dots\sigma_{k-2}\, u\,Y & \text{if $j=i$,}\\[0.5em]
D \sigma_j^{-1}\dots \sigma_1^{-1}\tau^a\sigma_1\dots\sigma_{k-2}\sigma_i\, u\,Y & \text{if $j>i$ and $j<k-1$,}\\[0.5em]
\mu_a\sigma_i\, u\,Y & \text{if $j=k-1$.}
\end{array}\right.
\end{equation}
It is now straightforward to check that $\Tr_k(X T^{(k)}_{j,a,u} Y)= X\Tr_k(T^{(k)}_{j,a,u})Y$ for $X=\sigma_i$, $i=1,\dots,k-2$.

A calculation similar to the one above and involving formulas (\ref{actiontau}) (with $\mathcal{V}_{j,\phi,u}$ replaced by $T^{(k)}_{j,a,u}$) shows that $\Tr_k(\tau^{\varepsilon} T^{(k)}_{j,a,u} Y)=\tau^{\varepsilon}\Tr_k(T^{(k)}_{j,a,u})Y$ where $\varepsilon=\pm 1$; we omit the details here and only indicate that 
one should consider separately the cases $j=0$, $0<j<k-1$ and $j=k-1$.

\vskip .2cm
(iii) Now we prove (\ref{cond-exp-val2b}). It is enough to consider the situation $k>1$, $X=T^{(k-1)}_{j,a,u}$, for $j=0,\dots,k-2$, $a\in \E_m$ and $u\in H(m,1,k-2)$. As $\sigma_{k-1}$ commutes with $u$, we have
\[\sigma_{k-1}^{-1}T^{(k-1)}_{j,a,u}\sigma_{k-1}=\sigma_{k-1}^{-1}\sigma_{j}^{-1}\dots\sigma_{1}^{-1}\tau^a\sigma_1\dots\sigma_{k-2}\sigma_{k-1}u=\left\{\begin{array}{ll}
T^{(k)}_{k-1,a,u} & \text{if $j=k-2$,}\\[0.5em]
T^{(k)}_{j,a,\sigma_{k-2}^{-1}u} & \text{if $j<k-2$.}
\end{array}\right.\]
Thus, with (\ref{def-tr}), we obtain
\[\Tr_k\left(\sigma_{k-1}^{-1}T^{(k-1)}_{j,a,u}\sigma_{k-1}\right)=\left\{\begin{array}{ll}
\mu_a\,u & \text{if $j=k-2$,}\\[0.5em]
D \sigma_j^{-1}\dots\sigma_1^{-1}\tau^a\sigma_1\dots\sigma_{k-3} u & \text{if $j<k-2$,}
\end{array}\right.\]
which is equal to $\Tr_{k-1}(T^{(k-1)}_{j,a,u})$ for any $j=0,\dots,k-2$.

Next, $\sigma_{k-1}X\sigma_{k-1}^{-1}=\sigma_{k-1}^{-1}X\sigma_{k-1}+(q-q^{-1})(X\sigma_{k-1} 
-\sigma_{k-1} 
X)$; by (\ref{cond-exp-val2a}) and (\ref{cond-exp-val2d}), 
$$\Tr_k(X\sigma_{k-1})=X\Tr_k(\sigma_{k-1})=DX=\Tr_k(\sigma_{k-1})X=\Tr_k(\sigma_{k-1}X)\ ,$$ 
for $X\in H(m,1,k-1)$, so $\Tr_k\left(\sigma_{k-1}X\sigma_{k-1}^{-1}\right)=\Tr_k\left(\sigma_{k-1}^{-1}X\sigma_{k-1}\right)=\Tr_{k-1}(X)$. 

\vskip .2cm
(iv)
It is left to prove (\ref{cond-exp-val2c}).
We use the following Lemma.
\begin{lemm}\hspace{-.2cm}.\hspace{.2cm}\label{lemm-tatb}
$\Tr_1\bigl(\Tr_2(\sigma_1^{-1}\tau^a\sigma_1\tau^b\sigma_1)\bigr)=D\mu_{a+b}\,$ 
for any integers $a$ and $b$.
\end{lemm}
\emph{Proof of the Lemma.} By induction on $c\geq0$ (the induction base, for $c=1$, is Lemma \ref{lem3}), one finds
\begin{equation}\label{tatc}
 \sigma_1^{-1}\tau^a\sigma_1\tau^c=\tau^c\sigma_1^{-1}\tau^a\sigma_1+(q-q^{-1})\sum\limits_{i=1}^c(\tau^{a+i}\sigma_1^{-1}\tau^{c-i}-\tau^i\sigma_1^{-1}\tau^{a+c-i})\,,\ \ \ \ \ \ c\geq0\,.
\end{equation}
Multiplying both sides by $\tau^{-c}$, we obtain
\begin{equation}\label{tat-c}
 \sigma_1^{-1}\tau^a\sigma_1\tau^{-c}=\tau^{-c}\sigma_1^{-1}\tau^a\sigma_1-(q-q^{-1})\sum\limits_{i=1}^c(\tau^{a+i-c}\sigma_1^{-1}\tau^{-i}-\tau^{i-c}\sigma_1^{-1}\tau^{a-i})\,,\ \ \ \ \ \ c\geq0\,.
\end{equation}
Multiplying (\ref{tatc}) and (\ref{tat-c}) by $\sigma_1$ from the right and using (\ref{cond-exp-val2a})-(\ref{cond-exp-val2b}), we find, for $c\geq0$,
\[\Tr_2(\sigma_1^{-1}\tau^a\sigma_1\tau^c\sigma_1)=\bigl(D-(q-q^{-1})\bigr)\tau^{a+c}+(q-q^{-1})\tau^c\mu_a+(q-q^{-1})\sum\limits_{i=1}^c(\tau^{a+i}\mu_{c-i}-\tau^i\mu_{a+c-i})\ ,\] 
\[\Tr_2(\sigma_1^{-1}\tau^a\sigma_1\tau^{-c}\sigma_1)=\bigl(D-(q-q^{-1})\bigr)\tau^{a-c}-(q-q^{-1})\tau^{-c}\mu_a+(q-q^{-1})\sum\limits_{i=1}^c(\tau^{a+i-c}\mu_{-i}-\tau^{i-c}\mu_{a-i})\ .\] 
Applying $\Tr_1$ to both sides of these two equalities, we find that $\Tr_1\bigl(\Tr_2(\sigma_1^{-1}\tau^a\sigma_1\tau^c\sigma_1)\bigr)=D\mu_{a+c}$ 
and $\Tr_1\bigl(\Tr_2(\sigma_1^{-1}\tau^a\sigma_1\tau^{-c}\sigma_1)\bigr)=D\mu_{a-c}\ $. 
\hfill$\square$

\vskip .2cm
Let $Z=T^{(k)}_{j,a,u}$ with $j\in\{0,\dots,k-1\}$, $a\in \E_m$ and $u\in H(m,1,k-1)$.
\begin{itemize}
\item Let first $j<k-1$. Note that $\sigma_{k-1}Z-Z\sigma_{k-1}=\sigma_{k-1}^{-1}Z-Z\sigma_{k-1}^{-1}$ and write
\[\sigma_{k-1}^{-1}Z-Z\sigma_{k-1}^{-1}=\sigma_{k-1}^{-1}\,\sigma_j^{-1}\dots\sigma_1^{-1}\tau^a\sigma_1\dots\sigma_{k-1}u-\sigma_j^{-1}\dots\sigma_1^{-1}\tau^a\sigma_1\dots\sigma_{k-1}u\sigma_{k-1}^{-1}\ .\] 
Using (\ref{cond-exp-val2a})-(\ref{cond-exp-val2b}), we find
\[\Tr_k(\sigma_{k-1}^{-1}Z-Z\sigma_{k-1}^{-1})=\Tr_{k-1}(\sigma_j^{-1}\dots\sigma_1^{-1}\tau^a\sigma_1\dots\sigma_{k-2})u-\sigma_j^{-1}\dots\sigma_1^{-1}\tau^a\sigma_1\dots\sigma_{k-2}\Tr_{k-1}(u)\ ,\]
which implies, by (\ref{cond-exp-val2a}),   
that $\Tr_{k-1}\bigl(\Tr_k(\sigma_{k-1}^{-1}Z-Z\sigma_{k-1}^{-1})\bigr)=0$, as required.

\vskip .2cm
 \item We treat now the situation $j=k-1$. 
Let $u=T^{(k-1)}_{l,b,w}$ with $l\in\{0,\dots,k-2\}$, $b\in\E_m$ 
and $w\in H(m,1,k-2)$. We have 
 \begin{equation}\label{pr-tr2}\Tr_k(\sigma_{k-1}Z)=D\sigma_{k-2}^{-1}\dots\sigma_1^{-1}\tau^a\sigma_1\dots\sigma_{k-2}\cdot \sigma_l^{-1}\dots\sigma_1^{-1}\tau^b\sigma_1\dots\sigma_{k-2}w\ .\end{equation}
 If $l<k-2$ we rewrite $\sigma_1\dots\sigma_{k-2}\cdot \sigma_l^{-1}\dots\sigma_1^{-1}$ as $\sigma_{l+1}^{-1}\dots\sigma_2^{-1}\cdot \sigma_1 
 \dots\sigma_{k-2}$, and obtain for the right hand side of (\ref{pr-tr2}) $D \sigma_{l}^{-1}\dots\sigma_1^{-1}\cdot \sigma_{k-2}^{-1}\dots\sigma_1^{-1}\tau^a\sigma_1\tau^b\sigma_1\dots\sigma_{k-2}\cdot \sigma_1\dots\sigma_{k-3}w$. We conclude that
 \begin{equation}\label{pr-tr3}\Tr_{k-1}\bigl(\Tr_k(\sigma_{k-1}Z)\bigr)=\left\{\begin{array}{ll}D\mu_{a+b}w & \text{if $l=k-2$,}\\[0.5em]
 D \sigma_l^{-1}\dots\sigma_1^{-1}\Tr_2(\sigma_1^{-1}\tau^a\sigma_1\tau^b\sigma_1)\sigma_1\dots\sigma_{k-3}w\hspace{0.3cm} & \text{if $l<k-2$.}\end{array}\right.\end{equation}
 Now we write
 $Z\sigma_{k-1}=\sigma_{k-1}^{-1}\dots\sigma_1^{-1}\tau^a\sigma_1\dots\sigma_{k-1}\cdot
 \sigma_l^{-1}\dots\sigma_1^{-1}\tau^b\sigma_1\dots\sigma_{k-2}\sigma_{k-1}w\,$.
This is equal to
 \[ \sigma_{l}^{-1}\dots\sigma_1^{-1}\cdot\sigma_{k-1}^{-1}\dots\sigma_1^{-1}\tau^a\sigma_1\tau^b\sigma_1\dots\sigma_{k-1}\cdot \sigma_1\dots\sigma_{k-2}w\ .\]
We use the already proved properties of $\Tr_k$ and find
\[\Tr_k(Z\sigma_{k-1})=\sigma_{l}^{-1}\dots\sigma_1^{-1}\Tr_2(\sigma_1^{-1}\tau^a\sigma_1\tau^b\sigma_1)\sigma_1\dots\sigma_{k-2}w\ .\]
We obtain finally that
\begin{equation}\label{pr-tr4}
\Tr_{k-1}\bigl(\Tr_k(Z\sigma_{k-1})\bigr)=\left\{\begin{array}{ll}\Tr_1\bigl(\Tr_2(\sigma_1^{-1}\tau^a\sigma_1\tau^b\sigma_1)\bigr)w & \text{if $l=k-2$,}\\[0.5em]
 D \sigma_l^{-1}\dots\sigma_1^{-1}\Tr_2(\sigma_1^{-1}\tau^a\sigma_1\tau^b\sigma_1)\sigma_1\dots\sigma_{k-3}w\hspace{0.3cm} & \text{if $l<k-2$.}\end{array}\right.\end{equation}
The comparison of (\ref{pr-tr3}) and (\ref{pr-tr4}), with the help of the Lemma \ref{lemm-tatb}, ends the verification for $j=k-1$.\hfill$\square$
\end{itemize}

\begin{rema}{\hspace{-0.2cm}.\hspace{0.2cm}}  \label{rema-condexp} 
{\rm Let  ${\mathfrak{B}}\subset {\mathfrak{A}}$ be a unital inclusion of associative unital algebras. 
Let $\Psi$ be a central form on ${\mathfrak{A}}$, non-degenerate on both ${\mathfrak{A}}$ and ${\mathfrak{B}}$. 
ÊRecall that the conditional expectation is the map $\epsilon\colon {\mathfrak{A}}\to {\mathfrak{B}}$ such that $\Psi(a,b)=\Psi(\epsilon(a)b)$ for all $a\in{\mathfrak{A}}$ and $b\in{\mathfrak{B}}$. Assume that 
$\Tr=\Tr_1\circ\dots\circ\Tr_{k-1}\circ\Tr_k$ is non-degenerate on $H(m,1,k)$ and $H(m,1,k-1)$. It then follows that the 
conditional expectation is the map $\Tr_k$. Indeed, to this end we have to show that 
\begin{equation}\label{reltrconexp}
\Tr(ZX)=\Tr(\Tr_k(Z)X)\ \ \ \ \text{for all}\ \ Z\in H(m,1,k)\ \text{and}\  X\in H(m,1,k-1)\ .\end{equation}
It follows from the following stronger fact, which is a direct consequence of (\ref{cond-exp-val2a}) and $\Tr_k(1)=1$,
$$\Tr_k(Z) X=\Tr_k(\Tr_k(Z) X)\ \ \ \ \text{for all}\ \ Z\in H(m,1,k)\ \text{and}\  X\in H(m,1,k-1)\ .$$
\hfill$\triangle$ 
}\end{rema}

\section{\hspace{-0.55cm}.\hspace{0.55cm}
$\mathcal{B}$-multiplicative central forms on the algebra $H(m,1,n)$}\label{sec-symm}

In this Section, we study the right regular action of the generators of $H(m,1,n)$ in the basis $\mathcal{B}$ 
introduced in Section \ref{sec-nf2}. We then use these formulas to define a family $L^{\gamma}_{\ n}$ of central forms on $H(m,1,n)$. 
Up to the standard involution of $H(m,1,n)$, these forms constitute a subset, corresponding to $D=0$, of the 
Markov traces discussed in Remark \ref{rema-mark}. 
The forms $L^{\gamma}_{\ n}$ have a certain ``multiplicativity" property with respect to the basis $\mathcal{B}$, 
the value of $L^{\gamma}_{\ n}$ on an element $T^{(n)}_{j,a,u}$ is a product of two factors, the first one depends on $j$ and $a$, the second one on $u$. For this reason, we call 
these forms {\it $\mathcal{B}$-multiplicative}. 
In general, for $D\neq 0$, this multiplicativity is lost. The multiplicativity gives additional means to work, and we establish the properties of the forms $L^{\gamma}_{\ n}$ separately and in a different, 
than in the previous section, 
manner. Then, using the fusion formula, from \cite{OPdA4}, for the algebras $H(m,1,n)$, 
we calculate, for finite $m$, the weights of these central forms.

\subsection{{\hspace{-0.50cm}.\hspace{0.50cm}}Right multiplication by the generators.} 

Recall that the elements $T^{(k)}_{j,a,u}\in H(m,1,k)$,
defined in (\ref{def-T}), form a basis of $H(m,1,k)$ when $u$ runs through a set of basis elements of $H(m,1,k-1)$,
$j=0,\dots,k-1$ and $a\in \E_m$. 
The basis $\mathcal{B}$ of $H(m,1,n)$ is obtained
recursively in this way,
starting from the basis $\{\tau^a\,|\  a\in \E_m\}$ of $H(m,1,1)$. 

Also recall that the formulas (\ref{actionsigma}) and (\ref{actiontau}), with $T^{(n)}_{j,a,u}$ instead of $\mathcal{V}_{j,\phi,u}$, give the left regular action of the generators of $H(m,1,n)$. 

The following Lemma gives the right regular action of the generators of $H(m,1,n)$ on elements 
$T^{(n)}_{j,a,u}$, spanning $H(m,1,n)$; we continue to use the notation $\beta_j$, see (\ref{def-betaj}).
\begin{lemm}{\hspace{-.2cm}.\hspace{.2cm}}\label{action-right}
Let $n>1$. We have
\begin{equation}\label{action-r}
T^{(n)}_{j,a,u}\cdot x=T^{(n)}_{j,a,ux}\,,\ \ \ \text{for any $x\in\{\tau,\tau^{-1},\sigma_1,\dots,\sigma_{n-2}\}$ and $u\in  H(m,1,n-1)$\,,}
\end{equation}
and, for elements $u=T^{(n-1)}_{k,b,w}$, $k\in\{0,\dots,n-2\}$, $b\in\mathbb{Z}$ and $w\in H(m,1,n-2)$, spanning $H(m,1,n-1)$,
\begin{equation}\label{action-r2}
\begin{array}{l} 
\hspace{-0.0cm}T^{(n)}_{j,a,u}\cdot\sigma_{n-1}=
(q-q^{-1})\Bigl(T^{(n)}_{j,a+b,\sigma_{k+1}\dots\sigma_{n-2}w}+S_a \Bigr)\\[1.5em]
\hspace{2cm}+\ \left\{\begin{array}{ll}
T^{(n)}_{k+1,b,\beta_{j+1}\tau^a\beta_{n-1}^{-1}w} & \text{if $j\leq k$,}\\[1em]
T^{(n)}_{k,b,\beta_{j}\tau^a\beta_{n-1}^{-1}w}-(q-q^{-1})T^{(n)}_{j,b,\beta_{k+1}\tau^a\beta_{n-1}^{-1}w}
 & \text{if $j> k$,}
\end{array}\right.
\end{array}
\end{equation}
where $$S_a:=\left\{\begin{array}{ll}
\sum\limits_{c=1}^a(T^{(n)}_{j,c,\beta_{k+1}\tau^{a+b-c}\beta_{n-1}^{-1}w}-T^{(n)}_{j,c+b,\beta_{k+1}\tau^{a-c}\beta_{n-1}^{-1}w})\quad & \text{if $a>0$\ ,}\\[1em]
0 & \text{if $a=0$\ ,}\\
\sum\limits_{c=1}^{-a}(T^{(n)}_{j,c+a+b,\beta_{k+1}\tau^{-c}\beta_{n-1}^{-1}w}-T^{(n)}_{j,c+a,\beta_{k+1}\tau^{b-c}\beta_{n-1}^{-1}w})\quad & \text{if $a<0$\ .}
\end{array}\right.$$
\end{lemm}
\emph{Proof.} The formula (\ref{action-r}) is immediate. As for (\ref{action-r2}), notice that $\sigma_{n-1}$ commutes with $w$, so
\[u\sigma_{n-1}=\sigma_k^{-1}\sigma_{k-1}^{-1}\dots\sigma_1^{-1}\tau^b\sigma_1\dots\sigma_{n-2}\sigma_{n-1}w=T^{(n)}_{k,b,w}\ ,\]
which gives
\[T^{(n)}_{j,a,u}\cdot\sigma_{n-1}=\sigma_j^{-1}\sigma_{j-1}^{-1}\dots\sigma_1^{-1}\tau^a\sigma_1\dots\sigma_{n-2}\sigma_{n-1}T^{(n)}_{k,b,w}\ .\]
If $a\geq0$, a straightforward calculation with the help of the formulas (\ref{actionsigma}), (\ref{actiontau}) (and (\ref{pui-tau})),
leads to (\ref{action-r2}). If $a<0$, we use the formula (implied by (\ref{pui-tau}))
$$ \tau^{a}\cdot T^{(n)}_{j,b,u}=(q-q^{-1}) \sum\limits_{c=1}^{-a}\Bigl(T^{(n)}_{0,c+a+b,\beta_j\tau^{-c}u}-T^{(n)}_{0,c+a,\beta_j\tau^{b-c}u}\Bigr)+T^{(n)}_{j,b,\tau^a  u}\ ,\ \ j>0\ ,$$
together with  (\ref{actionsigma}) and (\ref{actiontau}) to obtain (\ref{action-r2}).
\hfill$\square$

\subsection{{\hspace{-0.50cm}.\hspace{0.50cm}}$\mathcal{B}$-multiplicative central forms on $H(m,1,n)$.}

Fix a linear functional $\gamma$ with values in $\mathcal{A}_m$ on the space of polynomial functions in $\tau$
and set $\gamma_a:=\gamma(\tau^a)$ for $a\in\mathbb{Z}$. The linear functional $\gamma$ 
is determined by the values $\gamma_a$, $a\in\mathfrak{E}_m$.
We define a set of linear forms $L^{\gamma}_{\ n}\ :\ H(m,1,n)\to\mathcal{A}_m$, $n=0,1,\dots$,  
by the initial condition $L^{\gamma}_{\ 0}(1)=1$ and by the recursion, for $n>0$,
\begin{equation}\label{lin-form-rec}
L^{\gamma}_{\ n}(T^{(n)}_{j,a,u})=\delta_j^{
n-1}\gamma_a\, L^{\gamma}_{\ n-1}(u)\ \quad\text{for $j\in\{0,\dots,
n-1\}$, $a\in\mathfrak{E}_m$ and $u\in H(m,1,
n-1)$.}
\end{equation}
{}For the chain $H(m,1,0)\subset H(m,1,1)\subset\dots\subset H(m,1,n)\subset\dots$ of Hecke algebras, the embedding $H(m,1,k-1)\to H(m,1,k)$, adapted to the basis $\mathcal{B}$, 
is given by $H(m,1,k-1)\ni u\mapsto T^{(k)}_{k-1,0,u}$, so, if necessary, one can impose the normalization condition $\gamma_0=1$ and speak about a linear form $L^{\gamma}$ on the chain.  

\begin{prop}{\hspace{-.2cm}.\hspace{.2cm}}\label{sym-form}
{}{\rm (i)} For all $n=0,1,\dots$ the linear form $L^{\gamma}_{\ n}$ is central, that is
\begin{equation}\label{central-form}
L^{\gamma}_{\ n}(xx')=L^{\gamma}_{\ n}(x'x)\ \quad\text{for $x,x'\in H(m,1,n)$.}
\end{equation}
{\rm (ii)} The form $L^{\gamma}_{\ n}$ is invariant with respect to the anti-involution $\varpi$, $\varpi(xy)=\varpi(y)\varpi(x)$, of $H(m,1,n)$, identical on the generators  $\tau,\tau^{-1},\sigma_1,\dots,\sigma_{n-1}$. 
\end{prop}
\emph{Proof.} We prove the Proposition by induction on $n$. The result for $n=0$ and 1 is immediate ($H(m,1,1)$ is commutative). Let $n>1$. 

\vskip .2cm
(i) It is enough to prove (\ref{central-form}) for $x=\tau,\tau^{-1},\sigma_1,\dots,\sigma_{n-1}$ and
$x'=T^{(n)}_{j,a,u}$ with $j\in\{0,\dots,n-1\}$, $a\in\mathfrak{E}_m$ and $u\in H(m,1,n-1)$. 

Thus, let $x$ be one of the generators $\tau,\tau^{-1},\sigma_1,\dots,\sigma_{n-2}$. A direct analysis of the formula (\ref{actionsigma}), for $i=1,\dots,n-2$, and of the formula (\ref{actiontau}) (with $T^{(n)}_{j,a,u}$ instead of $\mathcal{V}_{j,\phi,u}$) leads to
\[ L^{\gamma}_{\ n}(xx')=\delta_j^{n-1}\gamma_a\, L^{\gamma}_{\ n-1}(xu)\ .\]
On the other hand, by (\ref{action-r})
we have $x'x=T^{(n)}_{j,a,ux}$ and thus
\[ L^{\gamma}_{\ n}(x'x)=\delta_j^{n-1}\gamma_a\, L^{\gamma}_{\ n-1}(ux)\ .\]
The formula (\ref{central-form}) follows by the induction hypothesis.

Now let $x=\sigma_{n-1}$.
The formula (\ref{actionsigma}) for $i=n-1$ yields
\[ L^{\gamma}_{\ n}(\sigma_{n-1}x')=\delta_j^{n-2}\gamma_a\, L^{\gamma}_{\ n-1}(u)\ .\]
Set $u=T^{(n-1)}_{k,b,w}$ where $k\in\{0,\dots,n-2\}$, $b\in\mathfrak{E}_m$ and $w\in H(m,1,n-2)$.
Using Lemma \ref{action-right}, formula (\ref{action-r2}), we obtain
\[\begin{array}{l}L^{\gamma}_{\ n}(x'\sigma_{n-1})=(q-q^{-1})\Bigl(\delta_j^{n-1}\delta_k^{n-2}\gamma_{a+b}\gamma_0\, L^{\gamma}_{\ n-2}(w)+L^{\gamma}_{\ n}(S_a)\Bigr)
\\[1em]
\hspace{2cm}+\ \left\{\begin{array}{ll}
\delta_{k}^{n-2}\delta_j^{n-2}\gamma_{b}\gamma_a\, L^{\gamma}_{\ n-2}(w) & \text{if $j\leq k$,}\\[1em]
-(q-q^{-1})\delta_{j}^{n-1}\delta_k^{n-2}\gamma_b\gamma_a\, L^{\gamma}_{\ n-2}(w)\quad & \text{if $j> k$.}\end{array}\right.
\end{array}\]
We have
$$L^{\gamma}_{\ n}(S_a)=\left\{\begin{array}{ll}
\sum\limits_{c=1}^a\delta_j^{n-1}\delta_k^{n-2}(\gamma_c\gamma_{a+b-c}-\gamma_{c+b}\gamma_{a-c})L^{\gamma}_{\ n-2}(w)\ \  & \text{if $a>0$\ ,}\\[1em]
0 & \text{if $a=0$\ ,}\\
\sum\limits_{c=1}^{-a}\delta_j^{n-1}\delta_k^{n-2}(\gamma_{a+b+c}\gamma_{-c}-\gamma_{a+c}\gamma_{b-c})L^{\gamma}_{\ n-2}(w)\ \  & \text{if $a<0$\ ,}
\end{array}\right.$$
which reduces, for any $a$, to $L^{\gamma}_{\ n}(S_a)=(\gamma_a\gamma_{b}-\gamma_{a+b}\gamma_0)\delta_j^{n-1}\delta_k^{n-2}L^{\gamma}_{\ n-2}(w)$. 
Thus
\[ L^{\gamma}_{\ n}(x'\sigma_{n-1})=\delta_{k}^{n-2}\delta_j^{n-2}\gamma_{b}\gamma_a\, L^{\gamma}_{\ n-2}(w),\]
which is equal to $\delta_j^{n-2}\gamma_a\, L^{\gamma}_{\ n-1}(u)$. This concludes the proof of the part (i).

\vskip .2cm
(ii) Let $y=T^{(n)}_{j,a,u}$. Using the centrality of $L^{\gamma}_{\ n}$, we write
\begin{equation}\label{leii}L^{\gamma}_{\ n}(\varpi(y))=L^{\gamma}_{\ n}(\varpi(u)\sigma_{n-1}\dots\sigma_1\tau^a\sigma_1^{-1}\dots\sigma_j^{-1})=L^{\gamma}_{\ n}(\sigma_{n-1}\dots\sigma_1\tau^a\sigma_1^{-1}\dots\sigma_j^{-1}\varpi(u))\ ,\end{equation}
which vanishes for $j\neq n-1$: 
\begin{equation}\label{leii2}\text{if $z\in H(m,1,n-1)$ then $L^{\gamma}_{\ n}(\sigma_{n-1}z)=L^{\gamma}_{\ n}(T^{(n)}_{n-2,0,z})=0$}\ .\end{equation} 
Let $l\in\{1,\dots,n-1\}$; let $\psi_1$ be an element in the 
subalgebra of $H(m,1,n)$ generated by $\sigma_1,\dots,\sigma_{l-1}$  and  $\psi_2$ 
an element in the subalgebra of $H(m,1,n)$ generated by $\sigma_{l+1},\dots,\sigma_{n-1}$. Then, with $z\in H(m,1,n-1)$, 
$$\begin{array}{l}L^{\gamma}_{\ n}(\sigma_{n-1}\dots\sigma_1\tau^a\psi_1(\sigma_l-\sigma_l^{-1})\psi_2 z)=(q-q^{-1})L^{\gamma}_{\ n}(\sigma_{n-1}\dots\sigma_1\tau^a\psi_1\psi_2 z)\\[1em]
\hspace{1cm}=(q-q^{-1}) L^{\gamma}_{\ n}(\sigma_{n-1}\dots\sigma_1\tau^a\psi_2\psi_1 z)
=(q-q^{-1}) L^{\gamma}_{\ n}(\psi_2^{\downarrow 1}\sigma_{n-1}\dots\sigma_1\tau^a\psi_1 z)\\[1em]
\hspace{1cm}=(q-q^{-1}) L^{\gamma}_{\ n}(\sigma_{n-1}\dots\sigma_1\tau^a\psi_1 z\psi_2^{\downarrow 1})=0\ .
\end{array}$$
Here we denoted by $\psi_2^{\downarrow 1}$ the element obtained from $\psi_2$ by the shift of the index of the generators, $\sigma_{l+1}\mapsto\sigma_l,\dots,\sigma_{n-1}\mapsto\sigma_{n-2}$. Then we used the centrality of  $L^{\gamma}_{\ n}$ and (\ref{leii2}). 

Therefore, the string $\sigma_1^{-1}\dots\sigma_{n-1}^{-1}$ appearing (for $j=n-1$) 
 in the argument of $L^{\gamma}_{\ n}$ in the last expression in (\ref{leii}) can be replaced by $\sigma_1^{-1}\dots\sigma_{n-2}^{-1}\sigma_{n-1}$, then by  $\sigma_1^{-1}\dots\sigma_{n-3}^{-1}\sigma_{n-2}\sigma_{n-1}$, ...\,, by $\sigma_1\dots\sigma_{n-1}$. So, 
 $L^{\gamma}_{\ n}(\varpi(y))$ equals 
$$\delta_j^{n-1}L^{\gamma}_{\ n}(\sigma_{n-1}\dots\sigma_1\tau^a\beta_{n-1}^{-1}\varpi(u))=
\delta_j^{n-1}
L^{\gamma}_{\ n}(\sigma_{n-1}^{-1}\dots\sigma_1^{-1}\tau^a\beta_{n-1}^{-1}\varpi(u))=\delta_j^{n-1}\gamma_aL^{\gamma}_{\ n-1}(\varpi(u))\ ,$$
where we used the Lemma \ref{auxile}. The proof of (ii) follows by the induction hypothesis. \hfill$\square$

\vskip .1cm
A more, than (\ref{lin-form-rec}),  general Ansatz: $L_n(T^{(n)}_{j,a,u})=\delta_j^{n-1}\gamma_a^{(n)}\, L_{n-1}(u)$, with a 
functional $\gamma^{(n)}$ for each $n$, gives no essentially
new central forms: the centrality implies that $\gamma_a^{(n)}\gamma_b^{(n-1)}=\gamma_b^{(
n)}\gamma_a^{(n-1)}$ for all $a,b$. 

\vskip .1cm
In contrast to the anti-involution $\varpi$, the involution $\iota$, defined by (\ref{iota}),
does not leave invariant the central form $L^{\gamma}_{\ n}$. Let 
\begin{equation}\label{def-iL}\iota(L^{\gamma}_{\ n}):=\iota^{(0)}\circ L^{\gamma}_{\ n}\circ\iota\ \end{equation}
where $\iota^{(0)}$ is the restriction of $\iota$, see (\ref{iota}), 
 to the ring $\mathcal{A}_m$. Then, already on the example of $H(m,1,2)$, one sees that, in general, $\iota(L^{\gamma}_{\ n})$ is not equal to $L^{\tilde{\gamma}}_{\ n}$ for any $\tilde{\gamma}$.

\begin{rema}{\hspace{-0.2cm}.\hspace{0.2cm}}\label{markov} 
{\rm For any linear functional $\gamma$ such that $\gamma_0=1$, the central form $\iota(L^{\gamma}_{\ n})$ is 
one of the Markov traces on the algebra $H(m,1,n)$ constructed in \cite{L}, see Remark \ref{rema-mark}, 
since $\iota(L^{\gamma}_{\ n})\bigl(1\bigr)=1$ and, for $k=2,\dots,n$,
$$\iota(L^{\gamma}_{\ n})\bigl(\sigma_{k-1}x\bigr)=0\ \ \ \ \text{and}\ \ \ \ \iota(L^{\gamma}_{\ n})\bigl(\sigma_{k-1}\dots\sigma_1\tau^a\sigma_1^{-1}\dots\sigma_{k-1}^{-1}x\bigr)=\iota^{(0)}(\gamma_{-a})\,\iota(L^{\gamma}_{\ n})\bigl(x\bigr)\ ,$$
where $a\in\mathfrak{E}_m$ and $x\in H(m,1,k-1)$. 
\hfill$\triangle$ 
}\end{rema}

\begin{rema}{\hspace{-0.2cm}.\hspace{0.2cm}}\label{base38pred}
{\rm Let $\B_{k}^{+}$ 
be the set of elements $t^{+(k)}_{j,a}$, $j=0,\dots,k-1$ and $\alpha\in\mathfrak{E}_m$, defined by 
\begin{equation}\label{def-t+}
t^{+(k)}_{j,a}:=\left\{\begin{array}{l}\sigma_{j+1}\sigma_{j+2}\dots \sigma_{k-1}\ ,\ \textrm{if $a=0$\ ,}\\[.3cm]
\sigma_j\sigma_{j-1}\dots \sigma_1\tau^a\sigma_1\dots \sigma_{k-1}\ ,\ \textrm{if $a\in \mathfrak{E}_m\backslash\{0\}$\ ,}
\end{array}\right. \end{equation} 
and let $T^{+(k)}_{j,a,u}:=t^{+(k)}_{j,a}u$ where $u\in H(m,1,
k-1)$. Then $\B_{n}^{+}\dots \B_{1}^{+}$ is a basis of $H(m,1,n)$ related to the normal form (\ref{normalform2'}).
By the Lemma \ref{auxile}, the recursion for the linear form $L^{\gamma}_{\ n}$ in terms of elements $T^{+(n)}_{j,a,u}$ is the same as in terms of elements 
$T^{(n)}_{j,a,u}$, 
$$L^{\gamma}_{\ n}(T^{+(n)}_{j,a,u})=\delta_j^{n-1}\gamma_a\, L^{\gamma}_{\ n-1}(u)\ ,$$
showing the multiplicativity of the forms $L^{\gamma}_{\ n}$ with respect to the basis $\B_{n}^{+}\dots \B_{1}^{+}$.
\hfill$\triangle$ 
}\end{rema} 
 
\begin{rema}{\hspace{-0.2cm}.\hspace{0.2cm}}\label{base38}{\rm 
Let $\B_{\ k}^{\gamma}$ be the set of elements 
$$\sigma_l\sigma_{l+1}\dots\sigma_{k-1}\ ,\ l=1,2,\dots,k,\ \ \text{and}\ \ \beta_j(\tau^a-\gamma_a)\beta_{k-1}^{-1}\ ,\ j\in\{0,\dots,k-1\}\ ,\ a\in\mathfrak{E}_m\backslash\{0\}\ .$$
It follows from the Corollary \ref{normalform-fin} that $\B_{\ n}^{\gamma}\dots \B_{\ 1}^{\gamma}$ 
is a basis of $H(m,1,n)$. 
The unit element $1$ of $H(m,1,n)$ belongs to this basis
and values of the linear form $L^{\gamma}_{\ n}$ on the basis elements are
\begin{equation}\label{vozoprfo}L^{\gamma}_{\ n}(1)=\gamma_0^n\quad\text{and}\quad L^{\gamma}_{\ n}(x)=0\ \text{for any other $x\in 
\B_{\ n}^{\gamma}\dots \B_{\ 1}^{\gamma}
$.}\end{equation}
If $\gamma_0\neq 0$ one can rescale $\gamma$ to set $\gamma_0=1$; with this choice, the basis 
$\B_{\ n}^{\gamma}\dots \B_{\ 1}^{\gamma}$ is quasi-symmetric, in the terminology of \cite{MM}, with respect to the form $L^{\gamma}_{\ n}$. 

\vskip .2cm
Fix the functional $\gamma^{\circ}$ by $\gamma^{\circ}_0=1$ and $\gamma^{\circ}_a=0$ for $a\in\mathfrak{E}_m\backslash\{0\}$. Then 
$\B_{\ k}^{\gamma^{\circ}}=\B_{k}$ and the basis $\B_{n}\dots \B_{1}$ is quasi-symmetric with respect to $L^{\gamma^{\circ}}_{\ n}$.  
The basis  $\B_{n}^+\dots \B_{1}^+$ from the Remark \ref{base38pred} is as well quasi-symmetric with respect to $L^{\gamma^{\circ}}_{\ n}$.  
Remark \ref{reduced}, together with the fact that the inverses of the generators do not appear in (\ref{def-t+}) 
implies in particular that the central form $L^{\gamma}$, for this special choice $\gamma=\gamma^{\circ}$, coincides with the central form on $H(m,1,n)$ defined, for finite $m$, in \cite{Bre-M}.
\hfill$\triangle$ 
}\end{rema}

\begin{rema}{\hspace{-0.2cm}.\hspace{0.2cm}} 
{\rm Each of the central forms $L^{\gamma}_{\ n}$ and $\iota(L^{\gamma}_{\ n})$  
becomes the canonical symmetrizing form on the group algebra of $G(m,1,n)$ when $q$ is specialized to 1, $\gamma$ to $\gamma^{\circ}$, and, for finite $m$, $v_j$ is 
specialized to $\xi_j$, $j=1,\dots,m$, where $\{\xi_j\}_{j=1,\dots,m}$ is the set of all $m$-th roots of unity. 
This implies that, for finite $m$, the central forms $L^{\gamma}_{\ n}$ and $\iota(L^{\gamma}_{\ n})$,
extended to $\mathbb{C}(q,v_1,\dots,v_m)\otimes_{\mathcal{A}_m} H(m,1,n)$, 
are non-degenerate for a generic choice of $\gamma$.
\hfill$\triangle$ 
}\end{rema} 

\subsection{{\hspace{-0.50cm}.\hspace{0.50cm}}Weights of $\mathcal{B}$-multiplicative central forms and fusion formula.}\label{weightsfusion}

In this Subsection we assume that $m<\infty$. We explain how the ``weights" of the central forms $L^{\gamma}_{\ n}$ and $\iota(L^{\gamma}_{\ n})$ are directly deduced from the fusion formula, obtained in \cite{OPdA4}, for the algebra $H(m,1,n)$. 

\paragraph{Preliminaries on multi-partitions.}
A  partition of $n$ is a tuple $\lambda=(\lambda_1,\dots,\lambda_l)$ of positive integers such that $\lambda_1\geqslant\lambda_2\geqslant\dots\geqslant\lambda_l$ and $n=\lambda_1+\dots+\lambda_l=:|\lambda|$.
We identify partitions with their Young diagrams; the Young diagram of $\lambda$ is a left-justified array of rows of
nodes containing $\lambda_j$ nodes in the $j$-th row, $j=1,\dots,l$; the rows are numbered from top to bottom
.
\vskip .2cm
An $m$-partition $\blambda=(\lambda^{(1)},\dots,\lambda^{(m)})$ of $n$ is an $m$-tuple of partitions such that $n=|\lambda^{(1)}|+\dots+|\lambda^{(m)}|$. 
We regard an $m$-partition as a set of  ``$m$-nodes"; an $m$-node $\balpha$ is a pair $(\alpha,k)$ consisting of a usual node $\alpha$ and an integer $k=1,\dots,m$, indicating to which diagram in the $m$-tuple the node belongs.

For an $m$-node $\balpha=(\alpha,k)$ lying in the line $x$ and the column $y$ of the $k$-th diagram, we set $\pos(\balpha):=k$ and $cc(\balpha):=y-x$. Let $q,v_1,\dots,v_m$ be the parameters of the algebra $H(m,1,n)$. We also set $c(\balpha):=v_kq^{2(y-x)}$ and call it the {\it quantum content} of the $m$-node $\balpha$.

\vskip .2cm
Let $\blambda$ be an $m$-partition. For $j=1,\dots,m$, let $\mathfrak{l}_{\blambda,x,j}$ be the number of nodes in the line $x$ of the $j$-th diagram of $\blambda$, and $\mathfrak{c}_{\blambda,y,j}$ 
the number of nodes in the column $y$ of the $j$-th diagram of $\blambda$. For an $m$-node $\balpha$ of $\blambda$, we define, as in \cite{CJ,OPdA4},  \emph{generalized hook lengths} $h^{(j)}_{\blambda}(\balpha)$, $j=1,\dots,m$, by
\[h^{(j)}_{\blambda}(\balpha):=\mathfrak{l}_{\blambda,x,j}+\mathfrak{c}_{\blambda,y,k}-x-y+1\,,\] 
if $\balpha$ lies in the line $x$ and the column $y$ of the $k$-th diagram of $\blambda$ (in particular, $h^{(k)}_{\blambda}(\balpha)$ is the usual hook length of $\alpha$ in $\lambda^{(k)}$).

Finally, we define
\begin{equation}\label{m-croc}
\textsf{F}_{\blambda}
:=(q^{-1}\!-q)^n\prod_{\balpha\in \blambda}\Biggl(c(\balpha)\prod_{k=1}^m\,\frac{q^{-cc(\balpha)}}{v_{\pos(\balpha)}q^{-h^{(k)}_{\blambda}(\balpha)}-v_kq^{h^{(k)}_{\blambda}(\balpha)}}\Biggr)\,.
\end{equation}
The element $\textsf{F}_{\blambda}$ can also be written as
\begin{equation}\label{m-croc2}
\textsf{F}_{\blambda}=\prod_{\balpha\in \blambda}\Biggl(\frac{q^{cc(\balpha)}}{\left[h_{\blambda}(\balpha)\right]_q}\prod_{\textrm{\scriptsize{$\begin{array}{c}k=1,\dots,m\\k\neq\pos(\balpha)\end{array}$}}}\frac{q^{-cc(\balpha)}}{v_{\pos(\balpha)}q^{-h^{(k)}_{\blambda}(\balpha)}-v_kq^{h^{(k)}_{\blambda}(\balpha)}}\Biggr)\,,
\end{equation}
where $[j]_q:=q^{j-1}+q^{j-3}+...+q^{-j+1}$ for a non-negative integer $j$.

\vskip .2cm
Let $\blambda$ be an $m$-partition of $n$. A standard $m$-tableau of shape $\blambda$ is a filling of $m$-nodes of
$\blambda$ with numbers $1,\dots,n$ in such a way that numbers in nodes increase rightwards in the lines and 
downwards in the columns in every diagram.

\paragraph{Weights of central forms.}
Let $\mathcal{F}_m:=\mathbb{C}(q,v_1,\dots,v_m)$ be the field of fractions of $\mathcal{A}_m$. The algebra  $\mathcal{F}_mH(m,1,n):=\mathcal{F}_m\otimes_{\mathcal{A}_m} H(m,1,n)$ is split semi-simple and its irreducible representations are in bijection with the set of $m$-partitions of $n$ \cite{AK}. For an $m$-partition $\blambda$ of $n$, we denote by $\chi_{\blambda}$ the associated irreducible character of $\mathcal{F}_mH(m,1,n)$. 

We naturally extend the central form $L^{\gamma}_{\ n}$ to a central form (that we still denote by $L^{\gamma}_{\ n}$) on $\mathcal{F}_mH(m,1,n)$. Due to the centrality, we have
$$L^{\gamma}_{\ n}=\sum \wbg\,\chi_{\blambda}\ ,$$
where the sum is over the set of $m$-partitions $\blambda$ of $n$. The coefficients $\wbg\in\mathcal{F}_m$ 
are called ``weights" of the central form $L^{\gamma}_{\ n}$. Similarly, we have
$$\iota(L^{\gamma}_{\ n})=\sum \twbg\,\chi_{\blambda}\ .$$
Using the explicit realization of the irreducible representations of $\mathcal{F}_mH(m,1,n)$ \cite{AK,OPdA},  one can verify that $\iota^{(0)}\circ \chi_{\blambda}\circ\iota=\chi_{\blambda}$, for any $m$-partition $\blambda$ of $n$ ($\iota$ and $\iota^{(0)}$ are extended respectively to $\mathcal{F}_mH(m,1,n)$ and $\mathcal{F}_m$ here). Thus, we have
\begin{equation}\label{w-tw}\wbg=\iota^{(0)}(\twbg)\ \ \ \ \ \ \text{for any $m$-partition $\blambda$ of $n$.}
\end{equation}

The weights $\twbg$ have been calculated in \cite{GIM} (the weights $\wbgo$ have been calculated independently in \cite{M}, see also \cite{CJ}). We will present different formulas for the weights $\wbg$ and $\twbg$ relying on the fusion procedure for the algebra $H(m,1,n)$.

\paragraph{Fusion formula for $H(m,1,n)$.} We briefly recall the fusion formula for the algebra $H(m,1,n)$ (see \cite{OPdA4} for more details).

\vskip .2cm
Define, for $i=1,\dots,n-1$, the {\it Baxterized} elements, with {\it spectral parameters} $\alpha$ and $\beta$:
\begin{equation}\label{bax-s}
\sigma_i(\alpha,\beta):=\sigma_i+(q-q^{-1})\frac{\beta}{\alpha-\beta}\ .
\end{equation}
Let $a_0,a_1,\dots,a_{m-1}\in\mathcal{A}_m$ be defined by 
\[(\rho-v_1)(\rho-v_2)\dots(\rho-v_m)=\rho^m+a_{m-1}\rho^{m-1}+\dots+a_1\rho+a_0\ ,\]
where $\rho$ is an indeterminate, and set
$\mathfrak{a}_i(\rho):=\rho^{m-i}+\rho^{m-i-1}a_{m-1}+\dots+\rho\,a_{i+1}+a_i$, for $i=0,\dots,m$. 
We define the following rational function with values in $H(m,1,n)$:
\begin{equation}\label{bax-t}
\tau(\rho):=\tau^{m-1}+\mathfrak{a}_{m-1}(\rho)\tau^{m-2}+\dots+\mathfrak{a}_2(\rho)\tau+\mathfrak{a}_1(\rho)\ .
\end{equation}
Finally, let 
$$\phi_{k}(u_1,\dots,u_k):=\sigma_{k-1}(u_k,u_{k-1})\sigma_{k-2}(u_k,u_{k-2})\dots \sigma_{1}(u_k,u_1)\tau(u_k)\sigma_1^{-1}\dots \sigma_{k-2}^{-1}\sigma_{k-1}^{-1}\ ,$$
and define the following rational function with values in the algebra $H(m,1,n)$:
\begin{equation}\label{def-Psi1}
\Phi(u_1,\dots,u_{n}):=
\phi_{n}(u_1,\dots,u_{n-1},u_{n})\phi_{n-1}(u_1,\dots,u_{n-1})\dots\dots \phi_{1}(u_1)
\ .
\end{equation}

\vskip .2cm
Let $\blambda$ be an $m$-partition of $n$ and ${\mathcal{T}}$ a standard $m$-tableau of shape $\blambda$. For $i=1,\dots,n$, we denote by $c_i$ the quantum content of the node of ${\mathcal{T}}$ containing $i$. The main result in \cite{OPdA4} is that the element $E_{{\mathcal{T}}}$, defined by the following consecutive evaluations
\begin{equation}\label{fus}
E_{{\mathcal{T}}}:=\mathsf{F}_{\blambda}\Phi(u_1,\dots,u_{n})\Bigr\rvert_{u_1=c_1}\dots\Bigr\rvert_{u_{n-1}=c_{n-1}}\Bigr\rvert_{u_{n}=c_{n}}\ ,
\end{equation}
is a primitive idempotent of $H(m,1,n)$ associated to the irreducible representation of $H(m,1,n)$ corresponding to $\blambda$. In particular, we have that
\begin{equation}\label{L-E}\twbg=\iota(L^{\gamma}_{\ n})\bigl(E_{{\mathcal{T}}}\bigr)\ .
\end{equation}

\begin{prop}\hspace{-.2cm}.\hspace{.2cm}\label{prop-weights}
We have:
\begin{equation}\label{weights}
{\rm{(i)}}
\quad\twbg=\mathsf{F}_{\blambda}\prod_{i=1}^n\bigl(\iota^{(0)}(\gamma_{-(m-1)})+\mathfrak{a}_{m-1}(c_i)\iota^{(0)}(\gamma_{-(m-2)})+\dots+\mathfrak{a}_2(c_i)\iota^{(0)}(\gamma_{-1})+\mathfrak{a}_1(c_i)\iota^{(0)}(\gamma_0)\bigr)\,,
\end{equation}
\begin{equation}\label{weights2}
{\rm{(ii)}}\quad
\twbg=\mathsf{F}_{\blambda}\ \iota^{(0)}\!\Biggl(\prod_{i=1}^n\bigl(c_i^{-m+1}\gamma_0 -\frac{c_i^{-m+2}}{a_0}\sum_{\mu=1}^{m-1}\mathfrak{a}_{\mu+1}(c_i)\gamma_{\mu}\bigr)\Biggr)\,.
\hspace{5.8cm}
\end{equation}
\end{prop}
\emph{Proof.} (i) We calculate $\twbg$ using (\ref{L-E}) and (\ref{fus}). The formula (\ref{weights}) is then a direct consequence of the definition 
(\ref{lin-form-rec}) and (\ref{def-iL}).

(ii) It is straightforward to check that $\tau(\rho)(\rho-\tau)=(\rho-v_1)(\rho-v_2)\dots(\rho-v_m)$, using that 
$\mathfrak{a}
_m(\rho)=1$ and the recurrence relation $\mathfrak{a}_{\mu}(\rho)=\rho\mathfrak{a}_{\mu+1}(\rho)+a_{\mu}$, $\mu=0,\dots,m-1$ (we also note that $\mathfrak{a}_0(\rho)=(\rho-v_1)(\rho-v_2)\dots(\rho-v_m)\,$). Thus, we have 
$$\tau(\rho)=\frac{(\rho-v_1)(\rho-v_2)\dots(\rho-v_m)}{\rho-\tau}\ .$$
A similar calculation yields 
$$(\rho-\tau^{-1})\bigl(\rho^{m-1}-\frac{\rho^{m-2}}{a_0}\sum_{\mu=1}^{m-1}\mathfrak{a}_{\mu+1}(\rho^{-1})\tau^{\mu}\bigr)=(\rho-v_1^{-1})(\rho-v_2^{-1})\dots (\rho-v_m^{-1})\ ,$$
where we used the formula $\tau^{-1}=-\displaystyle\frac{1}{a_0}(\tau^{m-1}+a_{m-1}\tau^{m-2}+\dots+a_2\tau+a_1)$ together with the formula $\displaystyle\frac{\rho^m\mathfrak{a}_0(\rho^{-1})}{a_0}=(\rho-v_1^{-1})(\rho-v_2^{-1})\dots (\rho-v_m^{-1})$. We deduce that, for any $\rho\in\mathcal{A}_m$,
$$\iota\bigl(\tau(\rho)\bigr)=\frac{(\widetilde{\rho}-v_1^{-1})(\widetilde{\rho}-v_2^{-1})\dots(\widetilde{\rho}-v_m^{-1})}{\widetilde{\rho}-\tau^{-1}}=\widetilde{\rho}^{m-1}-\frac{\widetilde{\rho}^{m-2}}{a_0}\sum_{\mu=1}^{m-1}\mathfrak{a}_{\mu+1}(\widetilde{\rho}^{-1})\tau^{\mu}\ \ \ \ \text{where $\widetilde{\rho}:=\iota^{(0)}(\rho)$\,.}$$
Using now this formula for $\iota\bigl(\tau(\rho)\bigr)$, we calculate $\twbg$ using  (\ref{fus})--(\ref{L-E}), (\ref{lin-form-rec}) and (\ref{def-iL}) as in item (i). We obtain (\ref{weights2}), noting that $\iota^{(0)}(c_i)=c_i^{-1}$, for $i=1,\dots,n$. \hfill$\square$

\vskip .2cm
Using (\ref{w-tw}), together with the fact that $\iota^{(0)}(\mathsf{F}_{\blambda})=\mathsf{F}_{\blambda}\,(-a_0)^n\displaystyle\prod_{i=1}^nc_i^{m-2}$ (which is obtained by a direct inspection of
(\ref{m-croc2}), recalling that $a_0=(-1)^mv_1\dots v_m$), we find that (\ref{weights2}) implies
\begin{equation}\label{weights3}\wbg=\mathsf{F}_{\blambda}\ \prod_{i=1}^n\Bigl(-\frac{a_0\gamma_0}{c_i} +\sum_{\mu=1}^{m-1}\mathfrak{a}_{\mu+1}(c_i)\gamma_{\mu}\Bigr)\,. 
\end{equation}

The explicit formulas for the weights yield the following criterion for the central forms $L^{\gamma}_{\ n}$ and $\iota(L^{\gamma}_{\ n})$ to be 
non-degenerate on $\mathcal{F}_mH(m,1,n)$.
\begin{coro}{\hspace{-0.2cm}.\hspace{0.2cm}}
The central forms $L^{\gamma}_{\ n}$ and $\iota(L^{\gamma}_{\ n})$ are non-degenerate on $\mathcal{F}_mH(m,1,n)$ if and only if 
$$-\frac{a_0\gamma_0}{c} +\sum_{\mu=1}^{m-1}\mathfrak{a}_{\mu+1}(c)\gamma_{\mu}\neq 0\ \ \ \ \ \ \text{for any $c\in\{v_pq^{\pm 2i}\ |\ p=1,\dots,m,\ i=0,1,\dots,n-1\}$.}$$
\end{coro}

\begin{rema}{\hspace{-0.2cm}.\hspace{0.2cm}} \label{tableau}
{\rm To calculate weights we evaluated the form on an idempotent $E_{{\mathcal{T}}}$, corresponding to an $m$-tableau. However, by centrality, the result depends only on 
the shape of the $m$-tableau.
Calculations for $D\neq 0$ are more difficult, {\it cf} a 
calculation for the usual Hecke algebra in \cite{IO2} where  
an evaluation on a particular tableau was used.
\hfill$\triangle$ 
}\end{rema}

\begin{rema}{\hspace{-0.2cm}.\hspace{0.2cm}} \label{schur}
{\rm In particular, if $\gamma=\gamma^{\circ}$ as in Remark \ref{base38}, the formula (\ref{weights3}) becomes:
$$\wbgo=(-a_0)^nc_1^{-1}\dots c_n^{-1}\mathsf{F}_{\blambda}=(q-q^{-1})^n\prod_{\balpha\in \blambda}\prod_{k=1}^m\frac{q^{-cc(\balpha)}}{q^{h^{(k)}_{\blambda}(\balpha)}-v_k^{-1}v_{\pos(\balpha)}q^{-h^{(k)}_{\blambda}(\balpha)}}\ .$$
We thus recover the so-called ``cancellation-free" formula obtained in \cite{CJ} for the Schur elements (the inverses of the weights) associated to the 
non-degenerate central form $L^{\gamma^{\circ}}_{\ n}$. The formulas (\ref{weights}), (\ref{weights2}) and (\ref{weights3}) are 
generalizations of this formula for the central forms $L^{\gamma}_{\ n}$ and $\iota(L^{\gamma}_{\ n})$. 
\hfill$\triangle$ 
}\end{rema}


\begin{thebibliography}{99}\addcontentsline{toc}{section}{References}

\bibitem{AK} Ariki S. and Koike K., \emph{A Hecke algebra of $(\mathbb{Z}/r\mathbb{Z})\wr S_n$ and construction of its irreducible representations}, Adv. in Math. 106 (1994) 216--243.

\bibitem{Bre-M} Bremke K. and Malle G., \emph{Reduced words and a length function for $G(e,1,n)$}, Indag. Mathem. 8 (1997), 453--469. 

\bibitem{Bro-M} Brou\'e M. and Malle G., \emph{Zyklotomische Heckealgebren}, Asterisque 212 (1993) 119--189.

\bibitem{C2} Cherednik I. V., \emph{A new interpretation of Gelfand--Tzetlin bases}, Duke Math. J. 54 (1987) 563--577. 

\bibitem{CJ} Chlouveraki M. and Jacon N., \emph{Schur elements for the Ariki--Koike algebra and applications},
 J. Algebraic Combin. 35(2) (2012) 291--311. ArXiv:1105.5910

\bibitem{CT} Coxeter H. and Todd J., \emph{A practical method for enumerating cosets of a finite abstract group}, Proc. Edinburgh Math. Soc. 5 (1936) 26--34.

\bibitem{GIM} Geck M., Iancu L. and Malle G., \emph{Weights of Markov traces and generic degrees}, Indag. Mathem. 11 (2000), 379--397.

\bibitem{IK} Isaev A. P.  and Kirillov A. N., \emph{Bethe subalgebras in Hecke algebra and Gaudin models}, ArXiv:1302.6495v1 [math.QA] 

\bibitem{IO} Isaev A. P.  and Ogievetsky O. V., \emph{On Baxterized solutions of reflection equation and integrable chain models}, Nuclear Physics B 760 (2007) 167--183. ArXiv: math-ph/0510078 

\bibitem{IO2}  Isaev A. P.  and Ogievetsky O. V., \emph{On representations of Hecke algebras}, Czechoslovak
J. Phys. 55 (2005) 1433--1441.

\bibitem{Kl} Kleshchev A., \emph{Linear and projective representations of symmetric groups}, Cambridge University Press, Cambridge, 2005.

\bibitem{L} Lambropoulou S., \emph{Knot theory related to generalized and cyclotomic Hecke algebras of type B}, J. Knot Theory and Its 
Ramifications 8 (1999) 621--658. ArXiv: math/0405504 [math.GT]

\bibitem{MM} Malle G. and Mathas A., \emph{Symmetric cyclotomic Hecke algebras}, J. Algebra 205 (1998) 275--293.

\bibitem{M} Mathas A., \emph{Matrix units and generic degrees for the Ariki--Koike algebras}, J. Algebra 281 (2004) 695--730.

\bibitem{MTV} Mukhin E., Tarasov V. and Varchenko A., \emph{Bethe subalgebras of the group algebra of the symmetric group}, Transformation Groups 18 (2013) 767--801.

\bibitem{OPdA} Ogievetsky O. and Poulain d'Andecy L., \emph{On representations of cyclotomic Hecke algebras}, Mod. Phys. Lett. A 26 No. 11 (2011) 795--803. ArXiv: 1012.5844 [math-ph]

\bibitem{OPdA2} Ogievetsky O. and Poulain d'Andecy L., \emph{Cyclotomic Hecke algebras: Jucys--Murphy elements, representations, classical limit}, preprint (2011). hal.archives-ouvertes.fr: hal-0064782

\bibitem{OPdA3a} Ogievetsky O. and Poulain d'Andecy L., \emph{Fusion procedure for Coxeter groups of type B and complex reflection groups $G(m,1,n)$}, to appear in Proc. Amer. Math. Soc. 
ArXiv: 1111.6293

\bibitem{OPdA3} Ogievetsky O. and Poulain d'Andecy L., \emph{On representations of complex reflection groups G(m,1,n)}, Theor. Math. Phys. 174 (2013) 95--108. 
ArXiv: 1205.3459 [math.RT] 

\bibitem{OPdA4} Ogievetsky O. and Poulain d'Andecy L., \emph{Fusion formula for cyclotomic Hecke algebras}, (2013).\\
ArXiv:1301.4237

\bibitem{ST} Shephard G.C. and Todd J.A., \emph{Finite unitary reflection groups}, Canad. J. Math. 6 (1954) 274--304.

\end{thebibliography}
\end{document}